\documentclass[12pt]{article}
\usepackage{theorem}
\usepackage{amssymb}
\usepackage{graphicx}
\usepackage[mathscr]{eucal}
\usepackage{amsbsy}
\usepackage{amsmath}
\newtheorem{prop}{}[section]

{\theorembodyfont{\upshape} \newtheorem{rema}[prop]{}}
\newcommand{\boma}[1]{{\mbox{\boldmath $#1$} }}

\begin{document}
\newcommand{\uper}[1]{\stackrel{\barray{c} {~} \\ \mbox{\footnotesize{#1}}\farray}{\longrightarrow} }
\newcommand{\nop}[1]{ \|#1\|_{\piu} }
\newcommand{\no}[1]{ \|#1\| }
\newcommand{\nom}[1]{ \|#1\|_{\meno} }
\newcommand{\uu}[1]{e^{#1 \AA}}
\newcommand{\UD}[1]{e^{#1 \Delta}}
\newcommand{\bb}[1]{\mathbb{{#1}}}
\newcommand{\HO}[1]{\bb{H}^{{#1}}}
\newcommand{\Hz}[1]{\bb{H}^{{#1}}_{\zz}}
\newcommand{\Hs}[1]{\bb{H}^{{#1}}_{\ss}}
\newcommand{\Hg}[1]{\bb{H}^{{#1}}_{\gg}}
\newcommand{\HM}[1]{\bb{H}^{{#1}}_{\so}}
\newcommand{\CM}[1]{\bb{C}^{{#1}}_{\so}}
\newcommand{\hz}[1]{H^{{#1}}_{\zz}}
\def\crit{crit}
\def\rbnw{0.51}
\def\rbnwp{0.52}
\def\reybnw{7.84}
\def\reybnwp{8.01}
\def\rtg{2.8}
\def\rtgp{2.9}
\def\reytg{5.07}
\def\reytgp{5.27}
\def\rkm{0.61}
\def\rkmp{0.62}
\def\reykm{1.00}
\def\reykmp{1.02}
\def\oh{O(d,\interi)}
\def\ohz{\oh \ltimes \Td}
\def\ohperz{\oh \times \Td}
\def\uno{{\bf 1}}
\def\t{{\tt{T}}}
\def\Aa{\mathscr{A}}
\def\Ti{\mathscr{T}}
\def\Isot{{\mathcal H}}
\def\PIsot{{\mathcal H}^{-}}
\def\Isott{{\mathcal H}_r}
\def\PIsott{{\mathcal H}^{-}_r}
\def\Rey{Re}
\def\ug{u^G}
\def\gaamma{\widehat{\gamma}}
\def\Ddd{\widehat{\Dd}}
\def\Rrr{\widehat{\Rr}}
\def\Tcc{\widehat{T}_{\tt{c}}}
\def\eep{\widehat{\epsilon}}
\def\mum{\hat{\mu}}
\def\uz{u_{*}}
\def\R{{R}\,}
\def\ti{{\tt{t}}}
\def\ef{\psi}
\def\fun{\mathcal{F}}
\def\fun{{\tt f}}
\def\tvainf{\vspace{-0.4cm} \barray{ccc} \vspace{-0,1cm}{~}
\\ \vspace{-0.2cm} \longrightarrow \\ \vspace{-0.2cm} \scriptstyle{T \vain + \infty} \farray}
\def\De{F}
\def\er{\epsilon}
\def\erd{\er_0}
\def\Tn{T_{\star}}
\def\Tc{T_{\tt{c}}}
\def\Tb{T_{\tt{b}}}
\def\Tl{\mathscr{T}}
\def\Tm{T}
\def\Ta{T_{\tt{a}}}
\def\ua{u_{\tt{a}}}
\def\Tg{T_{G}}
\def\Tgg{T_{I}}
\def\Tw{T_{w}}
\def\Ts{T_{\Ss}}
\def\Tr{\Tl}
\def\Sp{\Ss'}
\def\Tsp{T_{\Sp}}
\def\vsm{\vspace{-0.1cm}\noindent}
\def\comple{\scriptscriptstyle{\complessi}}
\def\nume{0.407}
\def\numerob{0.00724}
\def\deln{7/10}
\def\delnn{\dd{7 \over 10}}
\def\e{c}
\def\p{p}
\def\z{z}
\def\symd{{\mathfrak S}_d}
\def\del{\omega}
\def\Del{\delta}
\def\Di{\Delta}
\def\Ss{{\mathscr{S}}}
\def\Ww{{\mathscr{W}}}
\def\mmu{\hat{\mu}}
\def\rot{\mbox{rot}\,}
\def\curl{\mbox{curl}\,}
\def\Mm{\mathscr M}
\def\XS{\boma{x}}
\def\TS{\boma{t}}
\def\Lam{\boma{\eta}}
\def\DS{\boma{\rho}}
\def\KS{\boma{k}}
\def\LS{\boma{\lambda}}
\def\PR{\boma{p}}
\def\VS{\boma{v}}
\def\ski{\! \! \! \! \! \! \! \! \! \! \! \! \! \!}
\def\h{L}
\def\EM{M}
\def\EMP{M'}
\def\Rr{{\mathscr{R}}}
\def\Zz{{\mathscr{Z}}}
\def\E{{\mathscr E}}
\def\FFf{\mathscr{F}}
\def\A{F}
\def\Xim{\Xi_{\meno}}
\def\Ximn{\Xi_{n-1}}
\def\lan{\lambda}
\def\om{\omega}
\def\Om{\Omega}
\def\Sim{\Sigm}
\def\Sip{\Delta \Sigm}
\def\Sigm{{\mathscr{S}}}
\def\Ki{{\mathscr{K}}}
\def\Hi{{\mathscr{H}}}
\def\zz{{\scriptscriptstyle{0}}}
\def\ss{{\scriptscriptstyle{\Sigma}}}
\def\gg{{\scriptscriptstyle{\Gamma}}}
\def\so{\ss \zz}
\def\Dv{\bb{\DD}'}
\def\Dz{\bb{\DD}'_{\zz}}
\def\Ds{\bb{\DD}'_{\ss}}
\def\Dsz{\bb{\DD}'_{\so}}
\def\Dg{\bb{\DD}'_{\gg}}
\def\Ls{\bb{L}^2_{\ss}}
\def\Lg{\bb{L}^2_{\gg}}
\def\bF{{\bb{V}}}
\def\Fz{\bF_{\zz}}
\def\Fs{\bF_\ss}
\def\Fg{\bF_\gg}
\def\Pre{P}
\def\UU{{\mathcal U}}
\def\fiapp{\phi}
\def\PU{P1}
\def\PD{P2}
\def\PT{P3}
\def\PQ{P4}
\def\PC{P5}
\def\PS{P6}
\def\Q{P6}
\def\X{Q2}
\def\Xp{Q3}
\def\Vi{V}
\def\bVi{\bb{V}}
\def\K{V}
\def\Ks{\bb{\K}_\ss}
\def\Kz{\bb{\K}_0}
\def\KM{\bb{\K}_{\, \so}}
\def\HGG{\bb{H}^\G}
\def\HG{\bb{H}^\G_{\so}}
\def\EG{{\mathfrak{P}}^{\G}}
\def\G{G}
\def\de{\delta}
\def\esp{\sigma}
\def\dd{\displaystyle}
\def\LP{\mathfrak{L}}
\def\dive{\mbox{div}}
\def\la{\langle}
\def\ra{\rangle}
\def\um{u_{\meno}}
\def\uv{\mu_{\meno}}
\def\Fp{ {\textbf F_{\piu}} }
\def\Ff{ {\textbf F} }
\def\Fm{ {\textbf F_{\meno}} }
\def\Eb{ {\textbf E} }
\def\piu{\scriptscriptstyle{+}}
\def\meno{\scriptscriptstyle{-}}
\def\omeno{\scriptscriptstyle{\ominus}}
\def\Tt{ {\mathscr T} }
\def\Xx{ {\textbf X} }
\def\Yy{ {\textbf Y} }
\def\Ee{ {\textbf E} }
\def\VP{{\mbox{\tt VP}}}
\def\CP{{\mbox{\tt CP}}}
\def\cp{$\CP(f_0, t_0)\,$}
\def\cop{$\CP(f_0)\,$}
\def\copn{$\CP_n(f_0)\,$}
\def\vp{$\VP(f_0, t_0)\,$}
\def\vop{$\VP(f_0)\,$}
\def\vopn{$\VP_n(f_0)\,$}
\def\vopdue{$\VP_2(f_0)\,$}
\def\leqs{\leqslant}
\def\geqs{\geqslant}
\def\mat{{\frak g}}
\def\tG{t_{\scriptscriptstyle{G}}}
\def\tN{t_{\scriptscriptstyle{N}}}
\def\TK{t_{\scriptscriptstyle{K}}}
\def\CK{C_{\scriptscriptstyle{K}}}
\def\CN{C_{\scriptscriptstyle{N}}}
\def\CG{C_{\scriptscriptstyle{G}}}
\def\CCG{{\mathscr{C}}_{\scriptscriptstyle{G}}}
\def\tf{{\tt f}}
\def\ta{{\tt a}}
\def\tc{{\tt c}}
\def\tF{{\tt R}}
\def\C{{\mathscr C}}
\def\P{{\mathscr P}}
\def\V{{\mathscr V}}
\def\TI{\tilde{I}}
\def\TJ{\tilde{J}}
\def\Lin{\mbox{Lin}}
\def\Hinfc{ H^{\infty}(\reali^d, \complessi) }
\def\Hnc{ H^{n}(\reali^d, \complessi) }
\def\Hmc{ H^{m}(\reali^d, \complessi) }
\def\Hac{ H^{a}(\reali^d, \complessi) }
\def\Dc{\DD(\reali^d, \complessi)}
\def\Dpc{\DD'(\reali^d, \complessi)}
\def\Sc{\SS(\reali^d, \complessi)}
\def\Spc{\SS'(\reali^d, \complessi)}
\def\Ldc{L^{2}(\reali^d, \complessi)}
\def\Lpc{L^{p}(\reali^d, \complessi)}
\def\Lqc{L^{q}(\reali^d, \complessi)}
\def\Lrc{L^{r}(\reali^d, \complessi)}
\def\Hinfr{ H^{\infty}(\reali^d, \reali) }
\def\Hnr{ H^{n}(\reali^d, \reali) }
\def\Hmr{ H^{m}(\reali^d, \reali) }
\def\Har{ H^{a}(\reali^d, \reali) }
\def\Dr{\DD(\reali^d, \reali)}
\def\Dpr{\DD'(\reali^d, \reali)}
\def\Sr{\SS(\reali^d, \reali)}
\def\Spr{\SS'(\reali^d, \reali)}
\def\Ldr{L^{2}(\reali^d, \reali)}
\def\Hinfk{ H^{\infty}(\reali^d, \KKK) }
\def\Hnk{ H^{n}(\reali^d, \KKK) }
\def\Hmk{ H^{m}(\reali^d, \KKK) }
\def\Hak{ H^{a}(\reali^d, \KKK) }
\def\Dk{\DD(\reali^d, \KKK)}
\def\Dpk{\DD'(\reali^d, \KKK)}
\def\Sk{\SS(\reali^d, \KKK)}
\def\Spk{\SS'(\reali^d, \KKK)}
\def\Ldk{L^{2}(\reali^d, \KKK)}
\def\Knb{K^{best}_n}
\def\sc{\cdot}
\def\k{\mbox{{\tt k}}}
\def\x{\mbox{{\tt x}}}
\def\g{ {\textbf g} }
\def\QQQ{ {\textbf Q} }
\def\AAA{ {\textbf A} }
\def\gr{\mbox{gr}}
\def\sgr{\mbox{sgr}}
\def\loc{\mbox{loc}}
\def\PZ{{\Lambda}}
\def\PZAL{\mbox{P}^{0}_\alpha}
\def\epsilona{\epsilon^{\scriptscriptstyle{<}}}
\def\epsilonb{\epsilon^{\scriptscriptstyle{>}}}
\def\lgraffa{ \mbox{\Large $\{$ } \hskip -0.2cm}
\def\rgraffa{ \mbox{\Large $\}$ } }
\def\restriction{\upharpoonright}
\def\M{{\scriptscriptstyle{M}}}
\def\m{m}
\def\Fre{Fr\'echet~}
\def\I{{\mathcal N}}
\def\ap{{\scriptscriptstyle{ap}}}
\def\fiap{\varphi_{\ap}}
\def\dfiap{{\dot \varphi}_{\ap}}
\def\DDD{ {\mathfrak D} }
\def\BBB{ {\textbf B} }
\def\EEE{ {\textbf E} }
\def\GGG{ {\textbf G} }
\def\TTT{ {\textbf T} }
\def\KKK{ {\textbf K} }
\def\HHH{ {\textbf K} }
\def\FFi{ {\bf \Phi} }
\def\GGam{ {\bf \Gamma} }
\def\sc{ {\scriptstyle{\bullet} }}
\def\a{a}
\def\ep{\epsilon}
\def\c{\kappa}
\def\parn{\par \noindent}
\def\teta{M}
\def\elle{L}
\def\ro{\rho}
\def\al{\alpha}
\def\si{\sigma}
\def\be{\beta}
\def\ga{\gamma}
\def\te{\vartheta}
\def\ch{\chi}
\def\et{\eta}
\def\complessi{{\bf C}}
\def\len{{\bf L}}
\def\reali{{\bf R}}
\def\interi{{\bf Z}}
\def\Z{{\bf Z}}
\def\naturali{{\bf N}}
\def\Sfe{ {\bf S} }
\def\To{ {\bf T} }
\def\Td{ {\To}^d }
\def\Tt{ {\To}^3 }
\def\Zd{ \interi^d }
\def\Zt{ \interi^3 }
\def\Zet{{\mathscr{Z}}}
\def\Ze{\Zet^d}
\def\T1{{\textbf To}^{1}}
\def\es{s}
\def\ee{{E}}
\def\FF{\mathcal F}
\def\FFu{ {\textbf F_{1}} }
\def\FFd{ {\textbf F_{2}} }
\def\GG{{\mathcal G} }
\def\EE{{\mathcal E}}
\def\KK{{\mathcal K}}
\def\PP{{\mathcal P}}
\def\PPP{{\mathscr P}}
\def\PN{{\mathcal P}}
\def\PPN{{\mathscr P}}
\def\QQ{{\mathcal Q}}
\def\J{J}
\def\Np{{\hat{N}}}
\def\Lp{{\hat{L}}}
\def\Jp{{\hat{J}}}
\def\Pp{{\hat{P}}}
\def\Pip{{\hat{\Pi}}}
\def\Vp{{\hat{V}}}
\def\Ep{{\hat{E}}}
\def\Gp{{\hat{G}}}
\def\Kp{{\hat{K}}}
\def\Ip{{\hat{I}}}
\def\Tp{{\hat{T}}}
\def\Mp{{\hat{M}}}
\def\La{\Lambda}
\def\Ga{\Gamma}
\def\Si{\Sigma}
\def\Upsi{\Upsilon}
\def\Gam{\Gamma}
\def\Gag{{\check{\Gamma}}}
\def\Lap{{\hat{\Lambda}}}
\def\Upsig{{\check{\Upsilon}}}
\def\Kg{{\check{K}}}
\def\ellp{{\hat{\ell}}}
\def\j{j}
\def\jp{{\hat{j}}}
\def\BB{{\mathcal B}}
\def\LL{{\mathcal L}}
\def\MM{{\mathcal U}}
\def\SS{{\mathcal S}}
\def\DD{D}
\def\Dd{{\mathcal D}}
\def\VV{{\mathcal V}}
\def\WW{{\mathcal W}}
\def\OO{{\mathcal O}}
\def\RR{{\mathcal R}}
\def\TT{{\mathcal T}}
\def\AA{{\mathcal A}}
\def\CC{{\mathcal C}}
\def\JJ{{\mathcal J}}
\def\NN{{\mathcal N}}
\def\HH{{\mathcal H}}
\def\XX{{\mathcal X}}
\def\XXX{{\mathscr X}}
\def\YY{{\mathcal Y}}
\def\ZZ{{\mathcal Z}}
\def\CC{{\mathcal C}}
\def\cir{{\scriptscriptstyle \circ}}
\def\circa{\thickapprox}
\def\vain{\rightarrow}
\def\salto{\vskip 0.2truecm \noindent}
\def\spazio{\vskip 0.5truecm \noindent}
\def\vs1{\vskip 1cm \noindent}
\def\fine{\hfill $\square$ \vskip 0.2cm \noindent}
\def\ffine{\hfill $\lozenge$ \vskip 0.2cm \noindent}
\newcommand{\rref}[1]{(\ref{#1})}
\def\beq{\begin{equation}}
\def\feq{\end{equation}}
\def\beqq{\begin{eqnarray}}
\def\feqq{\end{eqnarray}}
\def\barray{\begin{array}}
\def\farray{\end{array}}
\makeatletter \@addtoreset{equation}{section}
\renewcommand{\theequation}{\thesection.\arabic{equation}}
\makeatother
\begin{titlepage}
{~}
\vspace{-2.5cm}
\begin{center}
{\Large \textbf{Large order Reynolds expansions \\ for the Navier-Stokes equations}}
\end{center}
\vspace{0.1truecm}
\begin{center}
{\large
Carlo Morosi$\,{}^a$, Mario Pernici $\,{}^b$, Livio Pizzocchero$\,{}^c$
({\footnote{Corresponding author}})
} \\
\vspace{0.3truecm}
{\small{
${}^a$ Dipartimento di Matematica, Politecnico
di Milano,
\\ P.za L. da Vinci 32, I-20133 Milano, Italy \\
e--mail: carlo.morosi@polimi.it \\
\vspace{0.2truecm} ${}^b$ Istituto Nazionale di Fisica Nucleare, Sezione di Milano, \\
Via Celoria 16, I-20133 Milano, Italy \\
e--mail: mario.pernici@mi.infn.it
\\
\vspace{0.2truecm} ${}^c$ Dipartimento di Matematica, Universit\`a di Milano\\
Via C. Saldini 50, I-20133 Milano, Italy\\
and Istituto Nazionale di Fisica Nucleare, Sezione di Milano, Italy \\
e--mail: livio.pizzocchero@unimi.it
}}
\end{center}
\begin{center}
{{\small \textbf{Abstract}}}
\end{center}
\vskip -0.2cm
{\small{
We consider the Cauchy problem for the incompressible homogeneous Navier-Stokes (NS) equations
on a $d$-dimen\-sional torus (typically, with $d=3$), in the $C^\infty$
formulation described, e.g., in \cite{smo}.
In the cited work and in \cite{appeul}
it was shown how to obtain quantitative estimates on the exact solution  of
the NS Cauchy problem via the \textsl{a posteriori} analysis
of an approximate solution; such estimates concern
the interval of existence of the exact solution
and its distance from the approximate solution,
evaluated in terms of Sobolev norms.
In the present paper we consider an approximate solutions of the NS Cauchy problem having the form
$u^N(t) = \sum_{j=0}^N R^j u_j(t)$,
where $R$ is the ``mathematical'' Reynolds number
(the reciprocal of the kinematic viscosity) and the coefficients
$u_j(t)$ are determined stipulating that the NS equations
be satisfied up to an error $O(R^{N+1})$. This subject
was already treated in \cite{apprey},
where, as an application, the Reynolds expansion of
order $N=5$ in dimension $d=3$ was considered for the initial datum
of Behr-Ne$\check{\mbox{c}}$as-Wu (BNW). In the present
paper, these results are enriched regarding both
the theoretical analysis and the applications. Concerning
the theoretical aspect, we refine the approach of
\cite{apprey} using results from \cite{smo}; moreover,
we show how to take into account the symmetries
of the initial datum in building up the expansion.
Concerning the applicative aspect we consider two more ($d=3$) initial
data, namely, the vortices of Taylor-Green (TG) and Kida-Murakami (KM);
the Reynolds expansions for the BNW, TG and KM data
are performed symbolically via a Python program, attaining
orders between $N=12$ and $N=20$.
Our a posteriori analysis indicates, amongst else, that the solution
of the NS equations with anyone of the above three data
is global if $R$ is below an explicitly computed critical value. Admittedly,
our critical Reynolds numbers are below the ones
characterizing the turbulent regime; however
these bounds have a sound theoretical support, are fully
quantitative and improve previous results of global existence.}}
\par
\vspace{0.2cm} \noindent
\textbf{Keywords:} Navier-Stokes equations, existence and regularity theory, theoretical approximation.
\textbf{AMS 2000 Subject classifications:} 35Q30, 76D03, 76D05.
\end{titlepage}
\section{Introduction and preliminaries}
\label{intro}
\textbf{Navier-Stokes (NS) equations; Reynolds number.}
The NS equations for an incompressible homogeneous fluid with no external forces, periodic
boundary conditions and initial datum $u_{*}$ can be written as
\beq {\partial u \over \partial \ti}  = \nu \Delta u + \PPP(u,u)~, \qquad u(x,0) = u_{*}(x)~.
\label{eulnu} \feq
Here: $\nu \in (0,+\infty)$ is the kinematic viscosity;
$u= u(x, \ti)$ is the divergence free velocity field;
the space variables $x = (x_s)_{s=1,...,d}$ belong to the torus
$\Td := (\reali/2 \pi \interi)^d$;
$\Delta := \sum_{s=1}^d \partial_{s s}$ is the Laplacian. Furthermore,
$\PPP$ is the bilinear map defined as follows:
for all sufficiently regular velocity fields $v, w$ on $\Td$,
\beq \PPP(v, w) := - \LP(v \sc \partial w) \label{ppp} \feq
where $(v \sc \partial w)_r := \sum_{s=1}^d v_s \partial_s w_r$ ($r=1,...,d$)
and $\LP$ is the Leray projection onto the space of divergence free vector fields.
As in \cite{apprey}, let us define
\beq t := \nu \ti~, \qquad R := {1 \over \nu}~; \label{setting} \feq
then Eq.\,\rref{eulnu} takes the form
\beq {\partial u \over \partial t}  = \Delta u + R \, \PPP(u,u)~, \qquad u(x,0) = u_{*}(x)~, \label{eul} \feq
to which we systematically refer in the sequel. In the present framework,
it is natural to define the Reynolds number as
\beq \Rey := {V_{*} L_{*} \over \nu} = V_{*} L_{*} R~, \label{derey} \feq
where $V_{*}$ is a characteristic
velocity and $L_{*}$ a characteristic length; later on we will give precise definitions for $V_{*}$ and $L_{*}$
as quadratic means related to the initial datum, which fit well to
our Sobolev framework (see Eqs.\rref{vs}-\rref{defrey}).
In the sequel $R$ and $\Rey$ will be referred to as the ``mathematical'' and the
``physical'' Reynolds number, respectively.
\salto
\textbf{NS functional setting.} The functional setting proposed in \cite{appeul}
for Eq.\,\rref{eul} was mainly based on $L^2$-type Sobolev spaces of finite
order; in \cite{smo} the attention passed to Sobolev spaces of infinite
order, made of $C^\infty$ functions. Both references
are relevant for our present purposes, so it convenient to review a few
issues from each one. \par
Let us start from the space $D'(\Td, \reali^d) \equiv \Dv$ of $\reali^d$-valued
distributions on $\Td$. Each $v \in \Dv$ has a weakly convergent Fourier
expansion $v = \sum_{k \in \Zd} v_k e_k$ where $e_k(x) := e^{i k \sc x}$ and
the coefficients $v_k \in \complessi^d$ fulfil the relations
$\overline{v_k} = v_{-k}$ due to the reality of $v$. The Laplacian
and the associated semigroup act on the whole space $\Dv$ and possess the
Fourier representations
\beq (\Delta v)_k = - |k|^2 v_k~, \qquad (e^{t \Delta} v)_k = e^{- t |k|^2} v_k
\label{repesp} \feq
($v \in \Dv$, $t \in [0,+\infty)$, $k \in \Zd$).
In the sequel  we consider the spaces $L^p(\Td, \reali^d) \equiv
\mathbb{L}^p$ for $p \in [1,+\infty)$; we are mainly interested in the case
$p=2$. For any $n \in\reali$,
the $n$-th Sobolev space of divergence free, zero mean vector
fields on $\Td$ is
\beq \HM{n}(\Td) \equiv \HM{n} := \{ v  \in \Dv~|~\dive \, v = 0,~
\langle v \rangle = 0,~ \sqrt{-\Delta}^{\,n} v \in \mathbb{L}^2 \}   \feq
$$ = \{ v  \in \Dv~| k \sc v_k = 0~ \mbox{for all $k$},~v_0 = 0,
\sum_{k \in \Zd \setminus \{0 \}} |k|^{2 n} |v_k| < + \infty \} $$
(in the above $\langle v \rangle$ indicates the mean over $\Td$, that
equals $v_0$). This is a Hilbert space
with the inner product
and the norm
\beq \la v | w \ra_n := \la \sqrt{-\Delta}^n v | \sqrt{-\Delta}^n w \ra_{L^2}
= (2 \pi)^d \! \! \! \! \! \sum_{k \in \Zd \setminus \{0 \}} \! \! \! |k|^{2 n} \overline{v_k} \sc w_k , \feq
\beq \| v \|_n := \sqrt{\la v | v \ra_n} \feq
($a \sc b := \sum_{r=1}^d a_r b_r$ for all $a, b \in \complessi^d$).
We have $\HM{p} \hookrightarrow \HM{n}$ if $n, p \in \reali$ and $n \leqs p$,
where $\hookrightarrow$ indicates a continuous embedding.
\par
We consider as well the infinite order Sobolev space
\beq \HM{\infty}(\Td) \equiv \HM{\infty}:= \cap_{n \in \reali} \HM{n} =
\cap_{n \in \naturali} \HM{n}~; \feq
this is a \Fre space with the locally convex topology induced by the family of norms $(\|~\|_n)_{n \in \reali}$
or, equivalently, by the countable subfamily $(\| ~\|_n)_{n \in \naturali}$.
For $k \in \naturali \cup \{ \infty \}$, let us consider the space
\beq \CM{k}(\Td) \equiv \CM{k} := \{ v \in C^k(\Td,\reali^d)~|~~
\dive v = 0,~ \la v \ra = 0~\}~, \feq
which is a Banach space for $k < \infty$ and a \Fre space for $k=\infty$, when
equipped with the sup norms for all derivatives up to order $k$.
Let $h, k \in \naturali$, $n \in \reali$;
then
\beq \CM{h} \hookrightarrow \HM{n}~~\mbox{if $h \geqs n$}~,
\qquad \HM{n} \hookrightarrow \CM{k}~~ \mbox{if $n > k + d/2$}, \feq
where the second statement depends on the Sobolev Lemma; these facts imply
\beq \HM{\infty} = \CM{\infty} \feq
(which indicates the equality of the above vector spaces and of
their \Fre topologies).
\par
Let us pass to the map $\PPP(v, w) := - \LP(v \sc \partial w)$ of
Eq.\,\rref{ppp}. The following facts are known:
\begin{itemize}
\item[(i)] the map $\PPP$ is well defined and bilinear from
$\HM{0} \times \HM{1}$ to $\LP \mathbb{L}^1_{0}$ (the image under the Leray projection
of the $L^1$, zero mean vector fields). In terms of Fourier coefficients,
\beq \PPP(v, w)_k = - i \LP_k \sum_{h \in \Zd} [v_h \sc (k-h)] w_{k - h}~, \label{furp} \feq
where $\LP_k : \complessi^d \to \complessi^d$ is the projection onto the orthogonal
complement of $k$.
\item[(ii)] For each real $n > d/2$, $\PPP$ sends continuously
$\HM{n} \times \HM{n+1}$ to $\HM{n}$; so, there is a constant $K_{n d} \equiv K_n$ such that
\beq \| \PPP(v, w) \|_n \leqs K_n \| v \|_n \| w \|_{n+1}
\quad \mbox{for $v \in \HM{n}$, $w \in \HM{n+1}$}~. \label{basineq} \feq
By the arbitrariness of $n$, one infers that $\PPP$ sends
continuously $\HM{\infty} \times \HM{\infty}$ to $\HM{\infty}$.
As a generalization of \rref{basineq}, for all real $p, n$ such that $p \geqs n > d/2$ there is a
constant $K_{p n d} \equiv K_{p n}$ such that
\beq \| \PPP(v, w) \|_p \leqs {1 \over 2} K_{p n} ( \| v \|_p \| w \|_{n+1} + \| v \|_n \| w \|_{p+1})
\label{basineqa} \feq
$$ \mbox{for $v \in \HM{p}$, $w \in \HM{p+1}$~.} $$
\item[(iii)] For each real $n > d/2+1$, there is a constant $G_{n d} \equiv G_n$
such that
\beq |\la \PPP(v,w) | w \ra_n | \leqs G_n \| v \|_n
\| w \|^{2}_{n} \quad \mbox{for $v \in \HM{n}$, $w \in \HM{n+1}$}~;
\label{katineq} \feq
this is the famous Kato inequality, see \cite{Kato}.
More generally, for all real $p, n$ such that $p \geqs n > d/2+1$ there is constant
$G_{p n d} \equiv G_{p n}$ such that
\beq | \la \PPP(v, w) | w \ra_p | \leqs
{1 \over 2} G_{p n} (\| v \|_p \| w \|_n + \| v \|_n \| w \|_p)\| w \|_p \label{katineqa} \feq
$$ \mbox{for $v \in \HM{p}$, $w \in \HM{p+1}$}~. $$
\end{itemize}
Papers \cite{cog} \cite{cok} give explicit (but probably non optimal) espressions
for the constants $K_n, G_n$ in the inequalities \rref{basineq} \rref{katineq}; in
particular, these references show that one can take
\beq K_3 = 0.323~, \qquad G_3 = 0.438 \qquad \mbox{if $d=3$}~. \label{k3g3} \feq
Eqs.\,\rref{basineqa} and \rref{katineqa} are ``tame'' refinements
(in the Nash-Moser sense) of \rref{basineq} and \rref{katineq}, respectively,
to which they are reduced for $p=n$ with $K_n = K_{n n}$ and $G_{n} = G_{n n}$.
Some relations very similar to these tame inequalities
have been used in  \cite{BKM} \cite{Tem}
and, more recently, in \cite{RSS}.
Appendix A of \cite{smo}, anticipating
a more detailed analysis to be presented
in \cite{coga}, gives explicit formulas for $K_{p n}$ and $G_{p n}$
for arbitrary $p, n$.
\vfill \eject \noindent
\textbf{The NS Cauchy problem, in a $\boma{C^\infty}$ formulation.} Let us choose
\beq R \in (0,+\infty)~, \qquad \uz \in \HM{\infty}~; \label{in}\feq
the corresponding NS Cauchy problem is:
\beq \mbox{Find}~
u \in C^\infty([0,T), \HM{\infty}) \quad \mbox{such that~~}
 {d u \over d t} = \Delta u + R \, \PPP(u,u), \quad u(0) = \uz~.  \label{cau} \feq
It is known that
\rref{cau} has a unique maximal (i.e., not extendable)
solution $u$, whose domain $[0, T)$ depends in principle
on $R$ and $\uz$; this gives by restriction any other solution.
Some classical references related to this subject
are \cite{BKM} \cite{Gig} \cite{Kato} \cite{Kat2}
\cite{Koz} \cite{Koz2} \cite{Tem} (some of these works
consider mainly a finite Sobolev order, but there
are standard arguments for passing to the infinite order
case, reviewed e.g. in \cite{smo}).
\salto
\textbf{The physical Reynolds number in terms of Sobolev norms.}
Let us return to Eq.\,\rref{derey}, defining the physical Reynolds number $\Rey$
in terms of some characteristic velocity $V_{*}$ and length $L_{*}$.
In this work we intend $V_{*}$ to be the initial mean quadratic velocity:
\beq V_{*} :=  \sqrt{{1 \over (2 \pi)^d} \int_{\Td} |u_{*}|^2 d x} = {1 \over (2 \pi)^{d/2}}
\| \uz \|_{L^2}~. \label{vs} \feq
Moreover we define the characteristic length $L_{*}$ as a quadratic mean
of $2 \pi/|k|$ over the Fourier modes of $\uz$, in the following way:
\beq L_{*} := \sqrt{{\sum_{k \in \Zd \setminus \{0 \}} (2 \pi/|k|)^2 | u_{* k} |^2 \over
\sum_{k \in \Zd \setminus \{0 \}} | u_{* k} |^2}} = 2 \pi {\| \uz \|_{-1} \over \| \uz \|_{L^2}}~. \label{elles}
\feq
Thus
\beq \Rey = V_{*} L_{*} R = {\| \uz \|_{-1} \over (2 \pi)^{d/2-1}} R~; \label{defrey} \feq
this will be our standard throughout the paper (note that, differently from \rref{elles},
\rref{defrey} makes sense for $\uz = 0$ as well). \par
In the sequel, for better convenience, any numerical
estimate involving the Reynolds number will be given for
both parameters $R$ and $\Rey$.
\salto
\textbf{The Reynolds expansion and its a posteriori analysis.}
For any $N \in \{0,1,2,...\}$, one can build an approximate
solution of the Cauchy problem \rref{cau} of the form
\beq u^N(t) = \sum_{j=0}^N R^j u_j(t) \feq
where the functions $u_0,...,u_N : [0,+\infty) \to \HM{\infty}$ are determined
so that $u^N(0)= u_{*}$ and ${d u^N/d t} - \Delta u^N - R \, \PPP(u^N,u^N) =
O(R^{N+1})$. A detailed analysis of this approximation has been presented
in \cite{apprey} in a slightly different framework, based on
Sobolev spaces of finite order; in the next section this construction
will be proposed in a version based on $\HM{\infty}$, and integrated
with some results for initial data having nontrivial symmetries. In a few words:
\begin{itemize}
\item[(i)] One has a recursion rule to compute $u_0, u_1,...~$. \parn
\item[(ii)] Once $u^N$ has been determined, it is possible to set up for it
an \textsl{a posteriori} analysis. The essential step in this direction
requires to choose a real $n > d/2+ 1$ and fix the attention on the Sobolev norms
\beq \| u^N(t) \|_n~, \quad \| u^N(t) \|_{n+1}~,
\label{norms2} \feq
\beq \| \big({d u^N \over d t} - \Delta u^N
- R \, \PPP(u^N, u^N))(t)\|_n~ \label{norms1} \feq
which measure the ``growth'' and the ``differential error'' of $u^N$
at order $n$ or $n+1$; a more refined analysis can be performed
considering, in addition to the norms \rref{norms2} \rref{norms1},
their analogues of any order  $p > n$. It is
important to remark that all the above norms can be explicitly computed
(or bounded from above) using only the known functions $u_0,...,u_N$. \par
Now, one applies to $u^N$ the general method
of \cite{appeul} (inspired by \cite{Che} \cite{Rob}) and \cite{smo} to get estimates
on the maximal solution $u$ of the problem \rref{cau} via
a posteriori analysis of any approximate solution. In this approach,
using the norms \rref{norms2} \rref{norms1} or some functions of time
which bind them from above, one writes down the so-called \emph{control Cauchy problem}:
this consists of a first order ODE for an unknown function $\Rr_n : [0,\Tc) \vain \reali$,
supplemented with the initial condition $\Rr_n(0) = 0$
(see the forthcoming Proposition \ref{mainpr}).
Assume this problem to have
a solution $\Rr_n$, with a suitable domain $[0, \Tc)$; then $\Rr_n$ is nonnegative,
the maximal solution $u$ of \rref{cau} has a domain larger that $[0, \Tc)$, and
\beq \| u(t) - u^N(t) \|_n \leqs \Rr_n(t) \qquad \mbox{for $t \in [0,\Tc)$}~. \label{buap} \feq
In particular, $u$ is global if $\Tc = +\infty$.
The solution $\Rr_n$ of the control Cauchy problem is typically
found numerically, using any standard package for the integration of ODEs.
After computing $\Rr_n$, for any $p > n$ one can use the analogues of order $p$ of
the norms \rref{norms2} \rref{norms1} (or some upper bounds of them)
to determine explicitly a function $\Rr_p: [0,\Tc) \vain \reali$
(again nonnegative) such that
\beq \| u(t) - u^N(t) \|_p \leqs \Rr_p(t) \qquad \mbox{for $t \in [0,\Tc)$}~. \label{buapp} \feq
\end{itemize}
\textbf{Plan of the paper and main results.} In Section \ref{due}, combining
results from \cite{apprey} and \cite{smo} we present the
Reynolds expansion in $\HM{\infty}$ and its a posteriori analysis via the control Cauchy problem.
Moreover we show that the symmetries of the initial datum are inherited by each term
$u_j$ in the Reynolds expansion; this fact allows to reduce the effort
in the recursive computation of the $u_{j}$'s. \par
In Section \ref{tre} we
present the Reynolds expansions for the ($d=3$) initial data of
Behr-Ne$\check{\mbox{c}}$as-Wu (BNW), Taylor-Green (TG) and Kida-Murakami (KM);
the expansions have been performed up to the orders
$N=20$ (for BNW and TG) and $N=12$ (for KM) allowing, amongst else,
to infer the global existence of the NS equations
\vfill \eject \noindent
{~}
\vskip -2cm \noindent
for $R \leqs \rbnw$, $R \leqs \rtg$ and $R \leqs \rkm$, respectively;
in terms of the physical Reynolds number \rref{defrey},
we have global existence if
$\Rey \leqs \reybnw$, $\Rey \leqs \reytg$
and $\Rey \leqs \reykm$, respectively.
In all these cases, the Reynolds expansions have been
computed symbolically using Python programs written for
this purpose; the orders $N=20$ or $N=12$
are the largest ones allowed by the PC we have
used to run these programs (see Section \ref{due} for
more details).
In the case
of the BNW datum, the present computations improve the results
obtained in \cite{apprey} with an expansion up to order $N=5$,
computed symbolically via Mathematica (on this point, see also item (v) in the next paragraph). \par
In each one of the above three cases, the symmetries of the initial datum have been employed
to reduce the computational costs. These symmetries are described
in Appendix A; they are particularly relevant for the KM datum,
that in fact arose in the papers by Kida and Murakami as a result
of their investigation on the highly symmetric vector fields
on $\Tt$.
\par
Finally, in Section 4 we present some speculations on how to push
to higher Reynolds numbers the present results of global existence
for the BNW, TG and KM data.
\salto
\textbf{An assessment of the previous quantitative bounds.}
We are aware that our upper bounds on $\Rey$
for global existence, of order $10$ at most, are much below
the numerical values of $\Rey$ related to
turbulence: for example in the classical papers by Brachet \textsl{et al} \cite{Bra1}
and Kida-Murakami \cite{Kid} \cite{Kid2}, where turbulence
is analyzed numerically for the TG and KM initial data,
the order of magnitude of $\Rey$ is between $10^2$ and $10^4$. However the cited
works, and all the other investigations on turbulence of which we are
aware, essentially assume global existence without proof. \par
We presume that our upper bounds on $\Rey$, being
supported by a rigorous theoretical analysis, may have some interest; their small
values correspond to the current state of the art concerning
global existence of strong NS solutions and are, in any case,
an improvement with respect to the estimates
one could derive from previous quantitative approaches
to the problem.
For a better understanding of the last statement let us
mention some results from earlier papers, and/or their implications
in the case of the BNW datum: \parn
(i) For $d=3$, paper \cite{Rob}
states global existence for the NS Cauchy problem in $\HM{1}$ whenever
$R \leqs 0.00724/\| \uz \|_1$ (see also the related
works \cite{Rob3} \cite{Rob2}).
\parn
(ii) Paper \cite{accau} improves the $d=3$ bound of \cite{Rob}
to $R \leqs 0.407/\| \uz \|_1$ (again, for
any initial datum in $\HM{1}$); for the BNW, TG and KM data
this improved condition reads $R \leqs 0.00527$,
$R \leqs 0.0298$, $R \leqs 0.00899$ respectively, or
(using the definition \rref{defrey}) $\Rey \leqs 0.0811$, $\Rey \leqs 0.0541$,
$\Rey \leqs 0.0147$. \parn
(iii) Again for $d=3$, papers [21][24] ensure global existence for the NS Cauchy problem
\rref{cau} whenever $R \leqs 1/(G_3 \| \uz \|_3) = 1/(0.438 \| \uz \|_3)$. In the BNW, TG and KM and cases,
this general estimate reads $R \leqs 0.0147$, $R \leqs 0.0557$, $R \leqs 0.00458$, respectively,
or $\Rey \leqs 0.227$, $\Rey \leqs 0.101$, $\Rey \leqs 0.00752$.
\parn
(iv) Paper \cite{appeul} uses a Galerkin approximant
and its a posteriori analysis to infer global
existence for the BNW datum under the condition
$R \leqs 0.125$, or $\Rey \leqs 1.92$. \parn
(v) We have already recalled that \cite{apprey}
presents an $N=5$ Reynolds expansion for the BNW datum;
this gives global existence when $R \leqs 0.23$,
or $\Rey \leqs 3.53$. \par
The above results on the BNW, TG and KM cases
are always weaker (or much weaker) than the present
outcomes $R \leqs \rbnw$, $R \leqs \rtg$, $R \leqs \rkm$
or $\Rey \leqs \reybnw$, $\Rey \leqs \reytg$, $\Rey \leqs \reykm$.
As already mentioned, in Section \ref{qua} we will indicate
some strategies that might yield future improvements.
\par
For completeness, let us also mention
papers \cite{CheG} \cite{Kuk} \cite{Rau} (and references
therein); these present (quantitative or semiquantitative)
conditions for global existence
of strong NS solutions in three dimensions with periodic
boundary conditions, for suitable initial data that,
in our language, would correspond to large values of $\Rey$;
however these data have small periods (i.e., fast oscillations) in one space direction.
On the contrary the BNW, TG and KM
data considered in this work are not highly
oscillating in any direction.
\section{The Reynolds expansion: recursion rules,
a posteriori analysis and symmetries}
From now on $K_n$, $G_n$ are contants fulfilling the inequalities
\rref{basineq} \rref{katineq}; $K_p, G_p$ and $K_{p n}, G_{p n}$
are constants fulfilling Eqs.\,\rref{basineq}\rref{katineq} with $n$ replaced by $p$, and
Eqs.\,\rref{basineqa}\rref{katineqa}.
\salto
\label{due}
\textbf{The expansion and its a posteriori analysis.}
Let us consider the NS Cauchy problem \rref{cau} with $R \in (0,+\infty)$
and a datum $\uz \in \HM{\infty}$.
We choose an order $N \in \{0,1,2,...\}$ and consider a function of the form
\parn
\vbox{
\beq u^N : [0, +\infty) \vain \HM{\infty}~,
\qquad t \mapsto u^N(t) := \sum_{j=0}^N R^j u_j(t)~, \label{defun} \feq
$$ u_j \in C^\infty([0,+\infty), \HM{\infty}) \qquad \mbox{for $j=0,...,N$}; $$
}
the functions $u_j$ herein are to be determined.
We regard $u^N$ as an ``approximate
solution'' of the NS Cauchy problem.
\begin{prop}
\label{un}
\textbf{Proposition.}
(i) Let $u^N$ be as in \rref{defun}; then
\beq {d u^N \over d t} - \Delta u^N - R \PPP(u^N, u^N) \label{eun} \feq
$$ = \Big({d u_0 \over d t} - \Delta u_0\Big) +
\sum_{j=1}^{N} R^j \Big[ {d u^j \over d t} - \Delta u_j - \sum_{\ell=0}^{j-1} \PPP(u_\ell, u_{j-\ell-1}) \Big]
- \sum_{j=N+1}^{2 N+1} R^j \! \! \! \! \! \sum_{\ell=j-N-1}^N \PPP(u_\ell, u_{j - \ell-1}). $$
(ii) One can define recursively a
family of functions $u_j \in C^\infty([0,+\infty), \HM{\infty})$ prescribing
the following, for $t \in [0,+\infty)$:
\beq u_0(t) := e^{t \Delta} \uz~, \label{recurz} \feq
\beq u_{j}(t) := \sum_{\ell=0}^{j-1}
\int_{0}^t d s \, e^{(t-s) \Delta} \PPP(u_{\ell}(s), u_{j-\ell-1}(s))
\quad (j=1,...,N)~. \label{recur} \feq
With this choice we have $u_0(0) = \uz$, $u_j(0) = 0$ for $j=1,...,N$ and
the coefficients of $R^0, R^1,..., R^N$ in the right hand side of Eq.\,\rref{eun} vanish, so that
\beq u^N(0) = \uz~, \label{daerr} \feq
\beq {d u^N \over d t} - \Delta u^N - R \PPP(u^N, u^N)
= - \sum_{j=N+1}^{2 N+1} R^j  \sum_{\ell=j-N-1}^N \PPP(u_\ell, u_{j - \ell-1})~.
\label{eunn} \feq
(iii) For any real $n > d/2$, Eqs.\,\rref{eunn}\rref{basineq} imply
the following for $t \in [0,+\infty)$:
\par
\vbox{
\beq \| \big({d u^N \over d t} - \Delta u^N - R \, \PPP(u^N, u^N)\big)(t) \|_n
\label{est} \feq
$$ \leqs K_n \sum_{j=N +1}^{2 N +1} R^j
\sum_{\ell=j-N -1}^N \| u_\ell(t) \|_n \| u_{j - \ell -1}(t) \|_{n+1}~.$$
}
\noindent
More generally, for any real $p \geqs n > d/2$, Eqs.\,\rref{eunn}\rref{basineqa} imply
\par
\vbox{
\beq \| \big({d u^N \over d t} - \Delta u^N - R \, \PPP(u^N, u^N)\big)(t) \|_p
\label{esta} \feq
$$ \leqs {1 \over 2} K_{p n} \sum_{j=N +1}^{2 N +1} R^j
\sum_{\ell=j-N -1}^N ( \| u_\ell(t) \|_p \| u_{j - \ell -1}(t) \|_{n+1} +
\| u_\ell(t) \|_n \| u_{j - \ell -1}(t) \|_{p+1})~. $$
}
\noindent
\end{prop}
\textbf{Proof.} An elementary variation of the proof of Proposition 3.1 in \cite{apprey}
(this considers the same subject for an initial datum $\uz$ in a finite order
Sobolev space, so that the functions $u_j, u^N$ have less regularity; the
adaptation to the present $\HM{\infty}$ framework is straightforward). \fine
\begin{rema}
\label{larem}
\textbf{Remarks.}
(i) Eqs.\,\rref{daerr} \rref{eunn} indicate the following:
$u^N$ satisfies the initial condition of the NS Cauchy problem
\rref{cau}, and it fulfils the evolution equation in \rref{cau} up to an error
described explicitly by \rref{eunn}. \parn
(ii) The recursive computation of $u_0, u_1,...,u_N$
via Eqs.\,\rref{recurz} \rref{recur} can be performed in terms of Fourier
coefficients; one uses the representations \rref{repesp} and
\rref{furp} for $e^{t \Delta}$ and $\PPP$. Due to the structure
of the recursion relations, the Fourier coefficients
of $u_0,u_1,...$ contain functions of time of the form $B_{a, b}(t) := t^a e^{- b t}$
with $a, b \in \naturali$; as already mentioned
in \cite{apprey}, the related computations involve integrals of the
form
\beq {~} \hspace{-0.5cm} \int_{0}^t \hspace{-0.3cm} d s \,  e^{-|k|^2 (t-s)} B_{a, b}(s) =
\left\{ \barray{ll}
\hspace{-0.2cm} a! \Big(  \dd{B_{0, |k|^2}(t) \over (b - |k|^2)^{a+1}}
- \sum_{\ell=0}^a {B_{\ell, b}(t) \over  (b - |k|^2)^{a + 1 - \ell} \ell!} \Big)
& \mbox{\hspace{-0.4cm} if $b \neq |k|^2$}; \\
\hspace{-0.2cm} \dd{B_{a + 1, |k|^2}(t) \over a + 1}  & \mbox{\hspace{-0.4cm} if $b = |k|^2$.}
\farray \right. \feq
The calculation of $u_0, u_1,... $ via the above rules
is particularly simple if the initial datum $\uz$ is a Fourier
polynomial, i.e., if it has finitely
many nonzero Fourier coefficients. In this case all the iterates
$u_0,...,u_N$ are Fourier polynomials as well,
and each one of their coefficients is a sum
$\sum_{a, b} C_{a, b} B_{a, b}(t)$ with $(a, b)$ in a
finite subset of $\naturali \times \naturali$
and $C_{a, b} \in \complessi^d$.
\fine
\end{rema}
Keeping in mind the previous facts, one can treat $u^N$ using the general framework
of \cite{appeul} \cite{smo} for approximate solutions of the NS equations; this analysis
of the Reynolds expansion was performed in \cite{apprey}
at the level of finite order Sobolev spaces, and its
analogue in $\HM{\infty}$ is presented hereafter.
Due to the strict connection between the present considerations and the inequalities
(\ref{basineq})-(\ref{katineqa}), here and in the sequel we always use
Sobolev norms of order $> d/2+1$.
\begin{prop}
\label{mainpr}
\textbf{Proposition.} (i) With $R, \uz$ as in \rref{in} and $N \in \{0,1,2,...\}$,
let $u^N$ and $u_0,...,u_N$ be as in
Eqs.\,\rref{defun} \rref{recurz} \rref{recur}; moreover, choose
a real $n > d/2+1$. Let $\Dd_n,
\Dd_{n+1}, \ep_n \in C([0,+\infty), [0,+\infty))$ be
growth and error estimators for $u^N$ of Sobolev orders $n$ or $n+1$, in the following sense:
\beq \| u^N(t) \|_n \leqs \Dd_n(t)~, \qquad \| u^N(t) \|_{n+1} \leqs \Dd_{n+1}(t)~, \label{ex0} \feq
\beq \| \big({d u^N \over d t} - \Delta u^N - R \, \PPP(u^N, u^N)\big)(t) \|_n \leqs
\ep_n(t)
\label{exx} \feq
for all $t \in [0,+\infty)$. Moreover, assume there is a function $\Rr_n \in C^1([0,\Tc), \reali)$
solving the ``control Cauchy problem''
\beq {d \Rr_n \over d t} = - \Rr_n
+ R (G_n \Dd_n + K_n \Dd_{n+1}) \Rr_n + R \, G_n \Rr^2_n + \ep_n~, \qquad
\Rr_n(0) = 0~. \label{cont} \feq
If $u \in C^\infty([0, T), \HM{\infty})$ is the maximal
solution of the NS Cauchy problem \rref{cau}, one has
\beq T \geqs \Tc~, \label{tta} \feq
\beq \| u(t) - \ua(t) \|_n \leqs \Rr_n(t) \qquad \mbox{for $t \in [0,\Tc)$}
\label{furth} \feq
In particular, $u$ is global ($T=+\infty$) if the control Cauchy problem \rref{cont}
has a global solution ($\Tc = +\infty$). \parn
(ii) Consider a real $p > n$ and let $\Dd_p, \Dd_{p+1}, \ep_p \in C([0,+\infty), [0,+\infty))$ be growth
and error estimators of orders $p$ or $p+1$, fulfilling
inequalities of the form \rref{ex0} \rref{exx} with $n$ replaced by $p$. Let
$\Rr_p \in C([0,\Tc),\reali)$ be the solution
of the linear problem
\beq {d \Rr_p \over d t} = - \Rr_p
+ R (G_p \Dd_p + K_p \Dd_{p+1} + G_{p n} \Rr_n) \Rr_p + \ep_p~, \qquad
\Rr_p(0) =0~, \label{contp} \feq
which is given explicitly by
\beq \Rr_p(t) = e^{\dd{-t + R \Aa_p(t)}} \int_{0}^t d s \, e^{\dd{s - R \Aa_p(s)}} \ep_p(s)~\qquad
\mbox{for $t \in [0,\Tc)$}~, \label{rp} \feq
$$ \Aa_p(t) := \int_{0}^ t d s \,
\big(G_p \Dd_p(s) + K_p \Dd_{p+1}(s) + G_{p n} \Rr_n(s)\big)~. $$
Then
\beq \| u(t) - \ua(t) \|_p \leqs \Rr_p(t) \qquad \mbox{for
$t \in [0,\Tc)$} \label{furthp} \feq
(incidentally, note that \rref{furth} \rref{furthp} imply $\Rr_n(t), \Rr_p(t) \geqs 0$).
\end{prop}
\textbf{Proof.} Use items (i)\,(ii) in Proposition 4.4 of \cite{smo}, choosing
as an approximate solution for the NS Cauchy problem the function $u^N$; note that
the physical time of \cite{smo}, say $\ti$, is related to the present
time variable by $t = \nu \ti = \ti/R$.
({\footnote{The cited reference indicates that the conditions
for $\Rr_n, \Rr_p$ in \rref{cont} \rref{contp} could be generalized
replacing everywhere the equality sign $=$ with $\geqs$; these generalizations
are not relevant for our present purposes.}})  \fine
\begin{rema}
\label{rema1}
\textbf{Remarks.} (i) In a few words, the previous proposition indicates
how to obtain bounds on the interval of existence of $u$ and on its
distance from $u^N$ of order $n$, or $p> n$, via an a posteriori analysis of $u^N$. Note
that the estimators $\Dd_n, \Dd_{n+1}, \ep_n$ appearing in the
control Cauchy problem \rref{cont} can be constructed using only $u^N$
(or its coefficients $u_0,...,u_N$); the same can be said
of the estimators $\Dd_p, \Dd_{p+1}, \ep_p$ in \rref{contp}.  \parn
(ii) The simplest choices for the above estimators are the tautological ones:
if $m$ is anyone of the Sobolev orders $n,n+1,p,p+1$ mentioned before, take for
$\Dd_m$ and $\ep_m$ the $m$-th norms of $u^N = \sum_{j=0}^N R^j u_j$ and of
$d u^N/d t - \Delta u^N - R \, \PPP(u^N, u^N)$,
as given by \rref{eunn}. One could consider alternative
estimators, which are rougher but computable with a smaller effort; these
have the form
\beq \Dd_m(t) := \sum_{j=0}^{N} R^j \| u_j(t) \|_m \label{choi} \feq
\beq \ep_m(t) := K_m \sum_{j=N +1}^{2 N +1} R^j
\sum_{\ell=j-N -1}^N \| u_\ell(t) \|_m \| u_{j - \ell -1}(t) \|_{m+1} \label{eest}~. \feq
The prescription \rref{eest} is suggested by Eq.\,\rref{est}; an
obvious variation of it for $m=p$ is suggested by Eq.\,\rref{esta}. The choices
\rref{choi} \rref{eest} reduce the construction
of the estimators to calculating the norms $\| u_j(t) \|_m$
or $\| u_j(t) \|_{m+1}$, which is less expensive than computing exactly the norms of
$u^N$ and $d u^N/d t - \Delta u^N - R \, \PPP(u^N, u^N)$.
An intermediate alternative is to compute exactly the involved norms
up to some order $M \in \{1,..., N\}$ in $R$, and bind the reminders more roughly;
this yields, for example, the estimators
\beq \Dd_m(t) := \| \sum_{j=0}^M R^j u_j(t) \|_m + \sum_{j=M+1}^{N} R^j \| u_j(t) \|_m~, \label{interm} \feq
that will appear in the applications of the next section. \parn
(iii) Eqs.\,\rref{furth}\rref{furthp} entail some rather obvious bounds on the difference between
the Fourier coefficients of $u(t)$ and $u^N(t)$; for
example, as mentioned in \cite{apprey}, Eq.\,\rref{furth} implies
\beq (2 \pi)^{d/2} | u_k(t) - u^N_k(t) | \leqs {\Rr_n(t) \over |k|^n}
\qquad \mbox{for $k \in \Zd \setminus \{0 \}$ and $t \in [0,\Tc)$}~. \label{eqmodi} \feq
\end{rema}
\textbf{Using the symmetries of the initial datum.} The present paragraph
is inspired by a setting proposed in \cite{bnw}
to treat symmetries of the incompressible Euler equations, that we are presently
adapting to the NS case (in any space dimension $d$).
We consider the group
$\oh$, formed by the orthogonal $d \times d$ matrices with integer entries:
\beq \oh := \{ S \in Mat(d \times d, \interi)~|~S^\t S = \uno_d \}~. \feq
A $d \times d$ matrix $S$ belongs to $\oh$ if and only if
\parn
\vbox{
\beq S = \mbox{diag}(\ep_1,...,\ep_d)  \, Q(\sigma) \label{rep} \feq
$$ \ep_s \in \{\pm 1\},~ Q(\sigma)~
\mbox{the matrix of a permutation $\sigma : \{1,...,d\} \vain \{1,...,d\}$}; $$}
more precisely, $Q(\sigma)$ is the matrix such that $(Q(\sigma) c)_s = c_{\sigma(s)}$
for all $c \in \complessi^d$, $s\in \{1,...,d\}$. Incidentally, the representation
\rref{rep} implies that each element of $S$ takes
values in $\{ -1, 0, 1 \}$. Counting the choices
for the signs $\ep_s$ and for $\sigma$ in Eq.\,\rref{rep}, one concludes
that $\oh$ has $2^d \times d!$
elements.
In particular, $O(3, \interi)$ has 48 elements; this group
is often indicated with $O_h$, and referred to as the \textsl{octahedral group}.
\par
To go on let us consider the semidirect product $\ohz$, i.e., the
Cartesian product
$\ohperz$, viewed as a group with the composition law
\beq (S,a) (U,b) := (S U, a + S b) \quad \qquad (S, U  \in \oh~;~ a, b \in \Td)~. \label{prcirc} \feq
Each element $(S,a) \in \ohz$ induces a rototranslation
\beq \E(S,a) : \Td \vain \Td~, \qquad x \mapsto \E(S,a)(x) := S x + a~, \feq
and the mapping $(S,a) \mapsto \E(S,a)$ is a group homomorphism
between $\ohz$ and the group of diffeomorphisms of $\Td$ into itself.  \par
Given a vector field
$v \in \HM{\infty}$ and $(S,a) \in \ohz$,
we can construct the push-forward
of $v$ along the mapping $\E(S,a)$; this is a
vector field $\E_{*}(S,a) v \in \HM{\infty}$ given by
\beq \E_{*}(S,a) v : \Td \vain \reali^d~, \qquad x \mapsto
(\E_{*}(S,a) \, v)(x) = S v(S^\t(x - a))~, \label{push} \feq
and its Fourier coefficients are
\beq (\E_{*}(S,a) \, v)_k = e^{-i a \sc k} S v_{\scriptscriptstyle{ S^\t k}}
 \qquad (k \in \Zd)~. \label{pushf} \feq
The linear map $\E_{*}(S,a) : v \mapsto \E_{*}(S,a) v$ preserves
the inner product $\la~|~\ra_n$ for each real $n$; denoting
with $O(\HM{\infty})$ the group of linear operators of $\HM{\infty}$
into itself preserving all inner products $\la~|~\ra_n$,
we have an injective group homomorphism
\beq \E_{*}: \ohz \vain O(\HM{\infty})~, \qquad (S, a) \mapsto \E_{*}(S,a)~. \feq
Let us turn the attention to the NS equations. Using \rref{pushf} with
the Fourier representations \rref{repesp} for $\Delta$,
$e^{t \Delta }$ and \rref{furp} for $\PPP$, one infers
\beq \Delta \E_{*}(S,a) \, v = \E_{*}(S,a) \Delta v~, \qquad
e^{t \Delta} \E_{*}(S,a) \, v = \E_{*}(S,a) e^{t \Delta} v~, \label{lape} \feq
\beq \PPP(\E_{*}(S,a) \, v, \E_{*}(S,a) \, w) = \E_{*}(S,a)\,  \PPP(v,w) \label{pinv} \feq
for $v,w \in \HM{\infty}$.
Now, consider the initial datum $\uz$ for the NS Cauchy problem \rref{cau} and
define the following, for $\sigma \in \{+,-\}$:
\beq \Isot^\si(\uz) := \{ (S,a) \in \ohz~|~\E_{*}(S,a) \uz = \si \uz \}~, \label{isot} \feq
\beq \Isott^\si(\uz) := \{ S \in \oh~|~(S,a) \in \Isot^\si(\uz)
~\mbox{for some $a \in \Td$} \}~. \label{isott} \feq
For $\sigma = +$, $\Isot^\si(\uz)$ and $\Isott^\si(\uz)$ are
subgroups of $\ohz$ and $\oh$; they will be
called the \emph{symmetry group} and the \emph{reduced symmetry group}
of $\uz$.
For $\sigma = -$, $\Isot^\si(\uz)$ and $\Isott^\si(\uz)$ will
be called the \emph{pseudo-symmetry} and \emph{reduced pseudo-symmetry
spaces} of $\uz$. The unions
$\Isot^{+}(\uz) \cup \Isot^{-}(\uz)$ and
$\Isott^{+}(\uz) \cup \Isott^{-}(\uz)$ are subgroups
of  $\ohz$ and $\oh$, respectively.
If $(\overline{S}, \overline{a})$ is any
element of $\Isot^{-}(\uz)$, then
$\Isot^{-}(\uz) = \Isot^{+}(\uz) \circ (\overline{S}, \overline{a})$
$ = (\overline{S}, \overline{a}) \circ \Isot^{+}(\uz)$
and $\Isott^{-}(\uz) = \Isott^{+} \, \overline{S}= \overline{S} \, \Isott^{+}(\uz)$
(here $\Isot^{+}(\uz) \circ (\overline{S}, \overline{a})$ means
$\{ (S,a) \circ (\overline{S}, \overline{a})~|~(S,a) \in \Isot^{+}(\uz) \}$, and so on).
\par
Using Eqs.\,\rref{lape} \rref{pinv}, one readily finds that
the iterates defined by
\rref{recurz} \rref{recur} fulfil at all times the relations
\beq \E_{*}(S,a) u_j(t) = \si^{j+1} u_j(t)\quad\mbox{for}~ (S, a)
\in \Isot^\si(\uz),
j \in \{0,1,..., N \}~. \feq
Recalling Eq.\,\rref{pushf}, one can rephrase
the above result in terms of Fourier coefficients, in the following way:
\beq u_{j, S k}(t) = \si^{j+1} e^{- i a \sc S k} S u_{j, k}(t)~
\mbox{for $(S,a) \in \Isot^\si(\uz)$, $j \in \{0,1,...,N\}$,
$k \in \Zd$}. \label{equal2} \feq
Let us point out the implications of the above results
in a concrete application of the recursion scheme \rref{recurz} \rref{recur},
say, for $d=3$, and with a Fourier polynomial as an initial datum;
in this case, Eq.\,\rref{equal2} can be employed to reduce the computational
cost for any iterate $u_j$. In fact, after computing
a Fourier coefficient $u_{j, k}$ one immediately obtains from
\rref{equal2} the coefficients $u_{j, S k}$ for all $S$ in
$\Isott^{+}(\uz) \cup \Isott^{-}(\uz)$: it suffices to apply the cited equation,
choosing for $a$ any element of $\Zd$ such that $(S,a) \in \Isot^{\pm}(\uz)$.
Note that $\{ S k \,|\, S \in \Isott^{+}(\uz) \cup \Isott^{-}(\uz)\}$ is
the orbit of $k$ with respect to the action of the group $\Isott^{+}(\uz) \cup \Isott^{-}(\uz)$
on $\Zt$.
\vfill \eject \noindent
\section{Applications. The BNW, TG and KM initial data}
\label{tre}
From here to the end of the paper, we consider the NS Cauchy problem
\rref{cau} with
\beq d = 3 ~; \feq
for the moment, the initial datum $\uz \in \HM{\infty}$ is unspecified. \par
We are interested in the Reynolds expansion $u^N(t) = \sum_{j=0}^N R^j u_j(t)$,
for suitable $N$, and on its a posteriori analysis via the control Cauchy problem
\rref{cont}. This will be performed choosing the Sobolev order
\beq n = 3~; \feq
thus Eq.\,\rref{cont} takes the form
\beq {d \Rr_3 \over d t} = - \Rr_3
+ R (G_3 \Dd_3 + K_3 \Dd_{4}) \Rr_3 + R \, G_3 \Rr^2_3 + \ep_3~, \qquad
\Rr_3(0) =0~, \label{conca} \feq
with $K_3$ and $G_3$ as in \rref{k3g3}; in the sequel we denote with
\beq \Rr_3 \in C^1([0, \Tc),\reali) \label{r3} \feq
the maximal solution, which is nonnegative.
\par
Let us recall that, after solving Eq. \rref{conca}, for any real $p > 3$
we could build a function $\Rr_p \in C^1([0,\Tc), \reali)$ following Eq. \rref{rp},
to be used in relation to estimates of Sobolev order $p$;
the actual computation of these higher order bounds,
in a number of applications, will be presented elsewhere.
\salto
\textbf{Choice of $\boma{\uz}$; automatic computations.}
As anticipated in the Introduction, in this paper we consider the BNW, TG and KM initial data;
these are Fourier polynomials,
described in detail in the sequel.
\par
For each one of these three data, the terms $u_0, u_1,..., $ in the Reynolds
expansion have been computed symbolically using Python
on a PC. To this purpose, we have developed the following software utilities:
\label{itemsab}
\begin{itemize}
\item[(a)] First of all, we have written
a Python program working in principle for any initial datum
$\uz$ of polynomial type; this computes
$u_0, u_1,...$ using Eqs.\,\rref{recurz} \rref{recur} and \rref{repesp} \rref{furp}.
\parn
\item[(b)] Secondly, for each one the BNW, TG and KM data
we have devised an \textsl{ad hoc} variant of the basic
program in (a), implementing the symmetries of the datum (see
Eq.\,\rref{equal2} and the related comments). These variants
reduce the computational costs, thus allowing to push
the Reynolds expansion to higher orders than the ones
allowed by the program in (a).
\end{itemize}
All the above Python programs use the package
GMPY \cite{gmpy} for fast arithmetics on rational numbers; they have been
run on an 8 Gb RAM PC.
Using the program mentioned in (a) we have attained
the orders $N = 16, 14, 7$ for the BNW, TG and KM datum, respectively.
Next we have used the specific Python programs mentioned in (b),
implementing the symmetries of these data; this has allowed
us to reach the orders $N=20,20,12$, respectively
({\footnote{As expected, type (b) programs give the same
result as the program of (a) up to the orders $16,14,7$;
this fact can be used to validate the implementation
of symmetries in these programs.
In the BNW case, it is also possible to make
a comparison with the $N=5$ expansion computed via Mathematica
in \cite{apprey}; again, there
is agreement between the results of the Mathematica and Python programs.}}).
The above Python programs also
give analytic expressions for the estimators $\Dd_3, \Dd_4$ and $\ep_3$
appearing in the control equation \rref{conca}.
For the three data mentioned before we have used
the estimator $\ep_3$ defined via \rref{eest}, and the
estimators $\Dd_3, \Dd_4$ defined via \rref{interm},
with $M=5$.
\par
The KM case has required the longest computational times;
calculations up to the order $N=12$ have taken, approximately,
$90$ hours for the determination of the $u_j$'s
and $30$ hours to compute the norms in Eqs.\,\rref{eest}
and \rref{interm}. Computations up to $N=20$ for the BNW and TG data
have been a bit faster, but in any case have required a few days.
\par
After computing the Reynolds expansion and the related
estimators, one can solve numerically the control Cauchy problem
\rref{conca}. This involves a Riccati type ODE with
time dependent coefficients, which have very long
analytic expressions when the order $N$ of
the expansion is large; for
the numerical treatment of this ODE
we have used Mathematica on a PC.
In our initial attempts, some numerical instabilities
have appeared for large $N$ in the integration of
\rref{conca}; these were due to insufficient precision
in the numerical evaluation of $\Dd_3, \Dd_4$ and
$\ep_3$ at the discrete times prescribed by the Mathematica
routines for ODEs. To eliminate
such instabilities, for the numerical
integration of \rref{conca} we have replaced $\Dd_3, \Dd_4$ and $\ep_3$
with convenient interpolants, built by
the internal routines of Mathematica after high
precision computations of the norms
in Eqs.\,\rref{eest} \rref{interm} at suitable grids of values for $t$.
For the choices of $N$ and $R$ considered
in our computations, the high precision computations
of the norms at a grid of instants and the
construction of the interpolants has required half an hour at most;
after this, the numerical solution of \rref{conca} has been almost instantaneous.
({\footnote{The results of the above numerical manipulations are
reliable but, admittedly, not certified.
Therefore, for the statements of the subsequent paragraphs stemming
from this approach to the ODE \rref{conca} we have a sound computer indication of validity, rather than a
``computer assisted proof''.
Perhaps, certified results could be obtained implementing
the calculation of $\Dd_3, \Dd_4, \ep_3$,
a finite-difference scheme for the ODE in \rref{conca}
and the related error analysis via any software for
interval arithmetics.}})
\salto
\textbf{General structure of the results from the control Cauchy problem.} Let
$\uz$ be the BNW, TG or KM datum. For any order $N$
considered in our computations, the numerical solution
of the control problem \rref{conca} for several values of
$R$ yields a picture already encountered in \cite{apprey}
for the BNW case and lower values of $N$. Recalling
that $[0, \Tc)$ is the domain of the maximal solution
$\Rr_3$ of \rref{conca}, we can summarize this picture
in the following way:
\begin{itemize}
\item[(i)] \label{pagitem} There is a critical number $R_{\crit}$, depending
on $\uz$ and $N$, such that
$\Tc = +\infty$ for $0 \leqs R \leqs R_{\crit}$, and $\Tc < + \infty$ for $R > R_{\crit}$.
Moreover, for $0 \leqs R \leqs R_{\crit}$ one has $\Rr_3(t) \vain 0^{+}$ for
$t \vain +\infty$, while for $R > R_{\crit}$ one has $\Rr_3(t) \vain + \infty$
for $t \vain \Tc^{-}$.
\item[(ii)] Let $u$ denote the maximal solution
of the NS Cauchy problem \rref{cau}. Due to (i),
for $0 \leqs R \leqs R_{\crit}$ $u$ is global and $\| u(t) - u^N(t) \|_3
\leqs \Rr_3(t)$ for all $t \in [0,+\infty)$; for $R > R_{\crit}$
it is only granted that $[0, \Tc)$ is in the domain of $u$, and
that $\| u(t) - u^N(t) \|_3 \leqs \Rr_3(t)$ for all $t \in [0,\Tc)$.
(Let us also recall Eq.\,\rref{eqmodi}, that can be used to bind
the Fourier coefficients $u(t) - u^N(t)$ via $\Rr_3(t)$.)
\item[(iii)] For the data and the values of $N$
considered in our computations (i.e., for $N$ up to a maximum
$20$ or $12$, depending on $\uz$),
$R_{\crit}$ increases with $N$.
\end{itemize}
One can associate to any $R$
a physical Reynolds number $\Rey$;
of course, item (ii) implies global existence for the NS Cauchy
problem when $\Rey$ is below the critical value $\Rey_{\crit}$,
defined as in \rref{defrey} with $R$ replaced by $R_{\crit}$.
For the largest values of $N$ attained in our
computations for the BNW, TG and KM data, $\Rey_{\crit}$
is close to the values anticipated in the Introduction, i.e.,
$\reybnw$, $\reytg$ and $\reykm$, respectively.
In the sequel we give more specific information analyzing
separately each one of the three initial data.
\vskip 0.2cm \noindent
\textbf{The BNW datum.} This is
\beq \uz(x_1, x_2, x_3) := 2 \big(\cos(x_1 + x_2) + \cos(x_1 + x_3), \label{unec0}\feq
$$ - \cos(x_1 + x_2) + \cos(x_2 + x_3),
 - \cos(x_1 + x_3) - \cos(x_2 + x_3)\big)~; $$
equivalently,
\parn
\vbox{
\beq \uz = \sum_{a=1}^3 z_a (e_{k_a} + e_{-{k_a}})~, \label{unec} \feq
$$ k_1 := (1,1,0),~~ k_2 := (1,0,1),~~ k_3 := (0,1,1)~; $$
$$ z_1 := (1,-1,0)~, \quad z_2 := (1,0,-1)~, \quad z_3 := (0,1,-1)~. $$
}
\parn
According to a conjecture of Behr, Ne$\check{\mbox{c}}$as and Wu
\cite{Nec}, this datum might produce a finite time blowup
for the Euler equations (i.e., for NS in the limit case of zero viscosity);
our position on this conjecture is described in \cite{bnw}. \par
Eqs.\,\rref{vs} \rref{elles} \rref{defrey} for this datum give $V_{*} = 2 \sqrt{3}$,
$L_{*} = 2 \pi/\sqrt{2}$ and
\beq \Rey = 2 \sqrt{6} \, \pi \, R = 15.39 ... \, R~. \feq
The symmetries of $\uz$ were already discussed in \cite{bnw}, and are
reviewed in Appendix A; in particular, the group $\Isott^{+}(\uz) \cup \Isott^{-}(\uz)$
has $12$ elements. We already mentioned the investigation of \cite{apprey}
on the Reynolds expansion for the BNW datum, performed
up to the order $N=5$ using Mathematica
({\footnote{The a posteriori analysis of \cite{apprey} was based on the tautological error estimators
$\Dd_m := \| u^5 \|_m$ ($m=3,4$) and
$\ep_3 := \| d u^5/d t - \Delta u^5 - R \PPP(u^5, u^5) \|_3$,
that could  be computed since the order $N=5$ is not too large. We repeat that,
for the higher order computations in the present work,
we have always used the rougher, but more easily computable
estimators in \rref{eest} with $M=5$.}}); the conclusion
of this analysis was a picture as in items (i)-(iii) before Eq.\rref{unec0}
where, for $N=5$, $R_{\crit} \in (0.23, 0.24)$ and, consequently, $\Rey_{\crit} \in (3.53, 3.70)$. \par
As already indicated, our Python program implementing
the BNW symmetries has allowed us to
push the expansion up to the order $N=20$. To give an idea of the computational
complexity we mention that $u_{20}$ has $6966$  nonzero Fourier coefficients,
whose wave vectors are partitioned
in $638$ orbits under the action of $\Isott^{+}(\uz) \cup \Isott^{-}(\uz)$ on $\Zt$.
Moreover, the nonzero Fourier coefficients of $u_{20}$ have very long expressions.
For example, let us consider $u^{(1)}_{20, k}$ for $k = (1,1,0)$, where $^{(1)}$ stands for the first of the three
components: $u^{(1)}_{20, (1,1,0)}(t)$ is a polynomial of degrees $9$ in $t$ and
$386$ in $e^{-t}$ with very complicated rational coefficients.
\par
Hereafter we summarize the results of the expansion up to $N=20$ and of its a
posteriori analysis via \rref{cont}. We have a picture as in the previously
cited items (i)-(iii) where, for $N=20$, \beq R_{\crit} \in (\rbnw, \rbnwp)
\quad  \mbox{whence} \quad \Rey_{\crit} \in (\reybnw, \reybnwp)~. \feq The
forthcoming Boxes 1a-1d present some results about computations with $N=20$ and
$R=\rbnw$, giving information on the following functions of time: the quantity
$(2 \pi)^{3/2} |u^{20}_k(t)|$ for the wave vector \hbox{$k=(1,1,0)$}; the
estimators $\Dd_3$ and $\ep_3$; the solution $\Rr_3$ of the control Cauchy
problem, which is global. In Boxes 2a-2d we consider the analogous functions in
the case $N=20$ and $R=\rbnwp$, in which $\Rr_3$ diverges at $\Tc = 2.855...$
({\footnote{An expression like $r= a. b c d e...~$ means that $a.b c d e$ are
the first digits of the output in the numerical computation of $r$.}}). Each
one of these boxes (and of the subsequent ones) contains the graph of the
function under consideration, and its numerical values for some choices of $t$.
\par Let us add a comment similar to one of \cite{apprey} about the pictures
that illustrated therein the $N=5$ BNW expansion. The functions in boxes of the
types (a) and (b) (i.e., $(2 \pi)^{3/2} |u^{20}_k(t)|$ and $\Dd_3(t)$) are very
similar in the cases $R = \rbnw$ and $R = \rbnwp$, even from the quantitative
viewpoint. Boxes (c) indicate that, as for $\ep_3$, the difference between the
cases $R = \rbnw$ and $R = \rbnwp$ is quantitatively significant; this is
sufficient to produce the completely different results for $\Rr_3$ illustrated
by boxes (d). Similar comments could be written about the boxes in the
forthcoming paragraphs, illustrating our computations about the TG and KM data.
\begin{figure}
\framebox{
\parbox{2in}{
\includegraphics[
height=1.3in,
width=2.0in
]%
{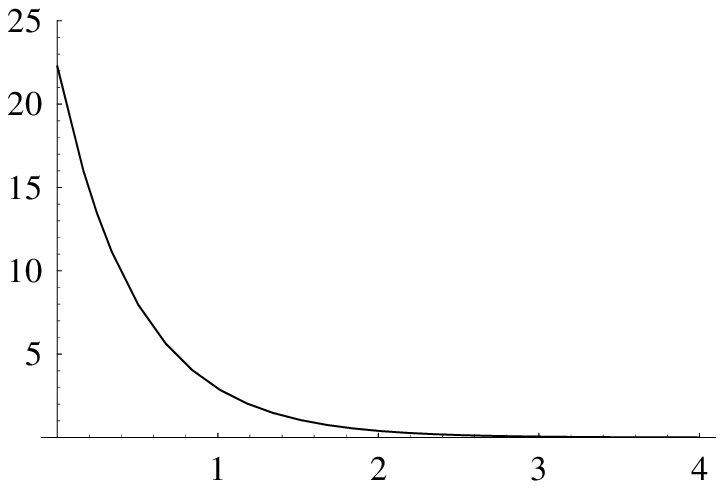}%
\par
{\tiny{ {\textbf{Box 1a.} BNW, $N=20$, $R = \rbnw$: the function
$\gamma(t) := (2 \pi)^{3/2} |u^{20}_{(1,1,0)}(t)|$. One has
$\gamma(0) = 22.27...$, $\gamma(0.5) = 8.031...$, $\gamma(1) =
2.933...$, $\gamma(1.5) = 1.077... $, $\gamma(2) = 0.396...$,
$\gamma(4) = 7.261... \times 10^{-3}$, $\gamma(8) = 2.435...
\times 10^{-6}$, $\gamma(10) = 4.461... \times 10^{-8}$~. }}
\par}
\label{f1a}
}
}
\hskip 0.4cm
\framebox{
\parbox{2in}{
\includegraphics[
height=1.3in,
width=2.0in
]%
{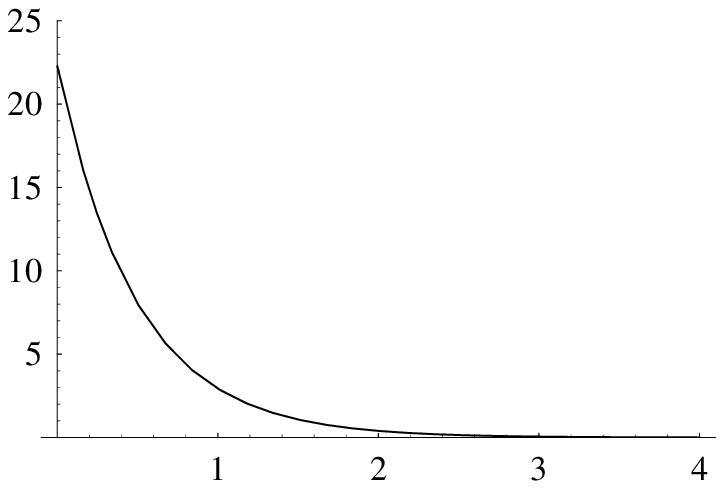}%
\par
{\tiny{ {\textbf{Box 2a.} BNW, $N=20$, $R=\rbnwp$: the function
$\gamma(t) := (2 \pi)^{3/2} |u^{20}_{(1,1,0)}(t)|$. One has
$\gamma(0) = 22.27...$, $\gamma(0.5) = 8.025...$, $\gamma(1) =
2.930...$, $\gamma(1.5) = 1.076... $, $\gamma(2) = 0.3960...$,
$\gamma(4) = 7.253... \times 10^{-3}$~. {~} \vskip 0.05cm
\noindent }}
\par}
\label{f2a}
}
}
\framebox{
\parbox{2in}{
\includegraphics[
height=1.3in,
width=2.0in
]%
{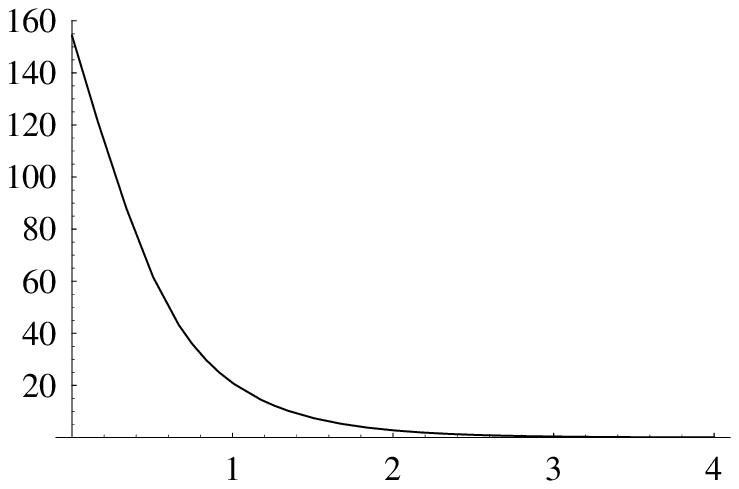}%
\par
{\tiny{ {\textbf{Box 1b.} BNW, $N=20, R=\rbnw$: the function
$\Dd_3(t)$. One has $\Dd_3(0) = 154.3...$, $\Dd_3(0.5) =
62.32...$, $\Dd_3(1) = 20.95...$, $\Dd_3(1.5) = 7.505... $,
$\Dd_3(2) = 2.748...$, $\Dd_3(4) = 0.05030...$, $\Dd_3(8)
=1.687... \times 10^{-5}$, $\Dd_3(10) = 3.091... \times 10^{-7}$~.
}}
\par}
\label{f1b}
}
}
\hskip 0.4cm
\framebox{
\parbox{2in}{
\includegraphics[
height=1.3in,
width=2.0in
]%
{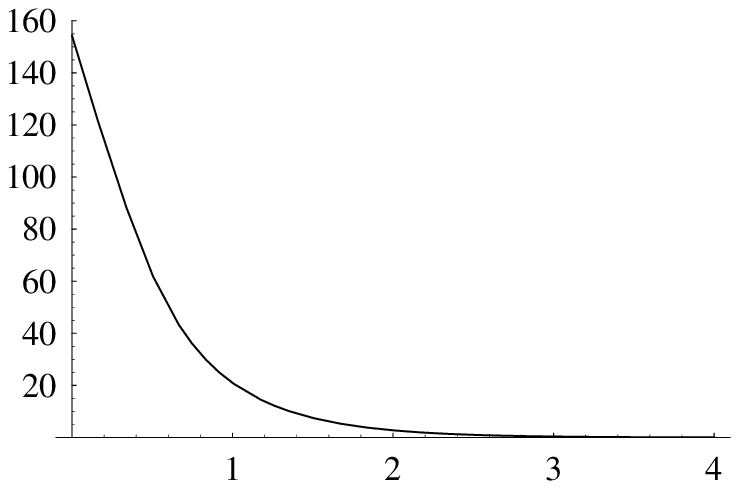}%
\par
{\tiny{
{\textbf{Box 2b.}
BNW, $N=20$, $R=\rbnwp$: the function
$\Dd_3(t)$.
One has $\Dd_3(0) = 154.3...$, $\Dd_3(0.5) = 62.53...$,
$\Dd_3(1) = 20.95...$, $\Dd_3(1.5) = 7.498... $, $\Dd_3(2) = 2.745...$,
$\Dd_3(4) = 0.05025...$~.
{~}
\vskip 0.2cm \noindent
}}
\par}
\label{f2b}
}
}
\framebox{
\parbox{2in}{
\includegraphics[
height=1.3in,
width=2.0in
]%
{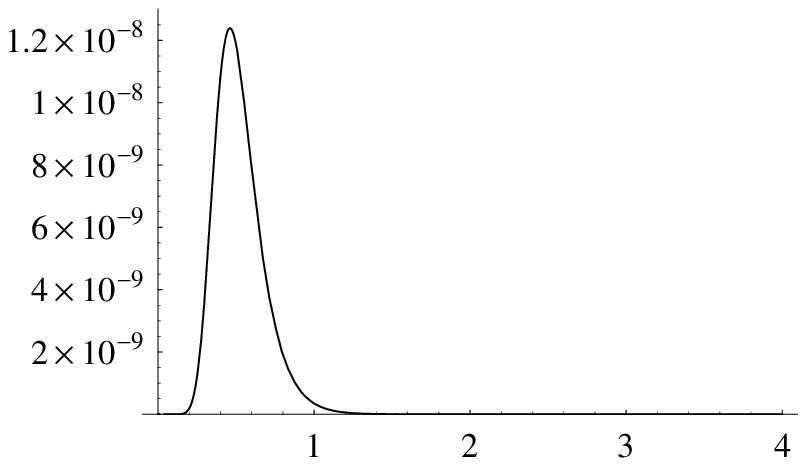}%
\par
{\tiny{
{\textbf{Box 1c.}
BNW, $N=20$, $R=\rbnw$: the function
$\ep_3(t)$.
One has $\ep_3(0) = 0$, $\ep_3(0.46) = 1.239 \times 10^{-8}...$,
$\ep_3(0.8) = 1.895 \times 10^{-9}...$,
$\ep_3(1) = 3.453 ...\times 10^{-10}$,
$\ep_3(2) = 5.753... \times 10^{-14}$, $\ep_3(4) = 5.695... \times 10^{-18}$\, .
}}
\par}
\label{f1c}
}
}
\hskip 0.4cm
\framebox{
\parbox{2in}{
\includegraphics[
height=1.3in,
width=2.0in
]%
{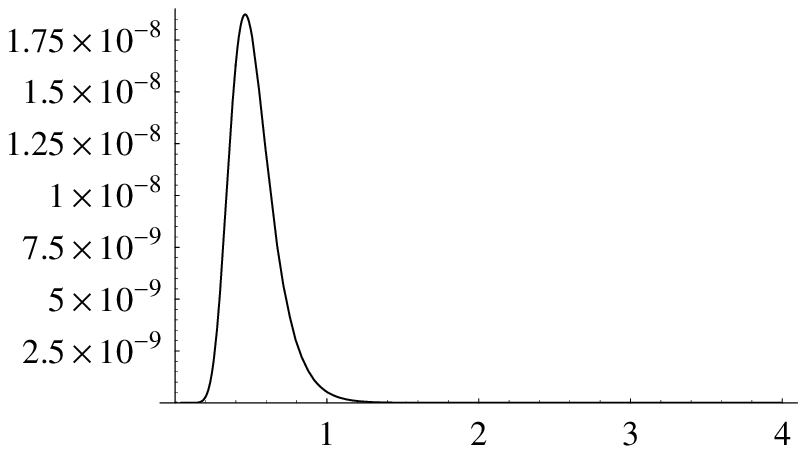}%
\par
{\tiny{
{\textbf{Box 2c.}
BNW, $N=20$, $R=\rbnwp$: the function $\ep_3(t)$.
One has $\ep_3(0) = 0$,
$\ep_3(0.46) = 1.868... \times 10^{-8}$,
$\ep_3(0.8) = 2.865...\times 10^{-9}$,
$\ep_3(1) = 5.219 ...\times 10^{-10}$,
$\ep_3(2) = 8.686... \times 10^{-14}$,
$\ep_3(4) = 8.572... \times 10^{-18}$\, .
\vskip 0.05cm
{~}
}}
\par}
\label{f2c}
}
}
\framebox{
\parbox{2in}{
\includegraphics[
height=1.3in,
width=2.0in
]%
{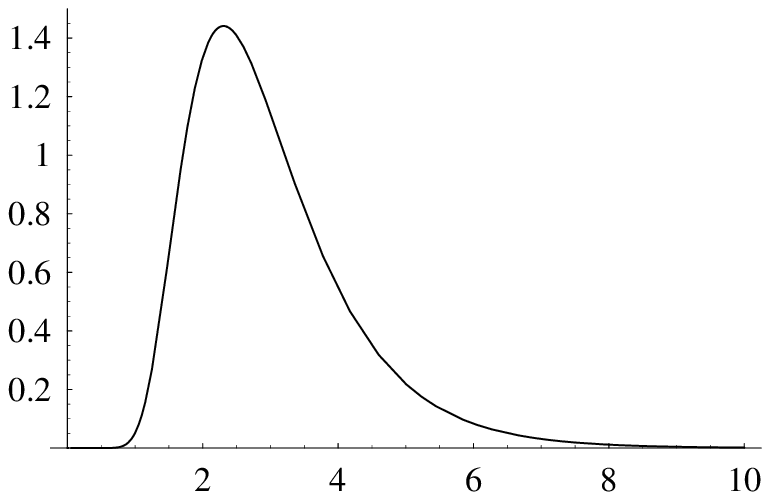}%
\par
{\tiny{
{\textbf{Box 1d.}
BNW, $N=20$, $R=\rbnw$: the function
$\Rr_3(t)$. This
appears to be \hbox{globally} defined, and vanishing at $+\infty$.
One has $\Rr_3(0) = 0$, $\Rr_3(1) = 0.05127...$,
$\Rr_3(1.5) = 0.6631...$, $\Rr_3(2.3) = 1.441...$, $\Rr_3(4) = 0.5433...$,
$\Rr_3(8) = 0.01143...$, $\Rr_3(10) = 1.551... \times 10^{-3}$~.
\vskip 0.1cm
{~}
}}
\par}
\label{f1d}
}
}
\hskip 3.7cm
\framebox{
\parbox{2in}{
\includegraphics[
height=1.3in,
width=2.0in
]%
{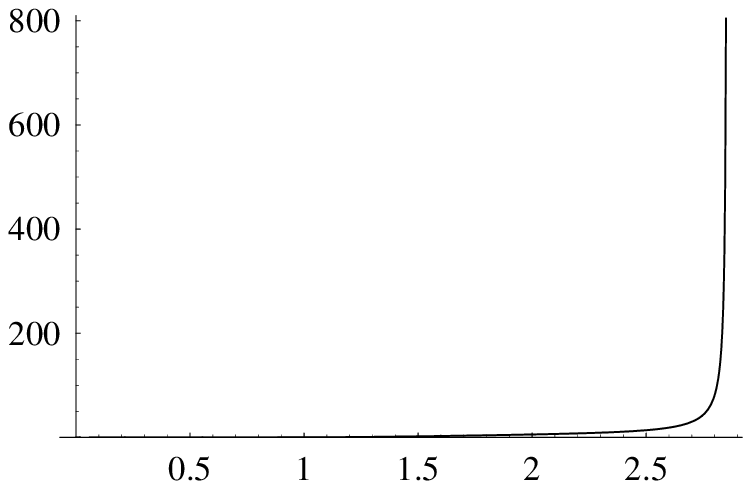}%
\par
{\tiny{
{\textbf{Box 2d.}
BNW, $N=20$, $R=\rbnwp$:
the function $\Rr_3(t)$. This diverges
for $t \vain \Tc = 2.855...$\,.
One has $\Rr_3(0) = 0$, $\Rr_3(0.5) = 1.847...\times 10^{-5}$,
$\Rr_3(1) = 0.1373...$, $\Rr_3(1.5) = 2.013...$,
$\Rr_3(2) = 5.611...$, $\Rr_3(2.85) = 804.5$~.
\vskip 0.35cm
{~}
}}
\par}
\label{f2d}
}
}
\end{figure}
\vfill\eject\noindent
\textbf{The TG datum.} This is \beq \uz(x_1, x_2, x_3) :=
\big(\sin x_1 \cos x_2 \cos x_3 ,  -\cos x_1 \sin x_2 \cos x_3 , 0 \big)~;
\label{tg} \feq equivalently, \beq \uz = {i \over 8} \sum_{a=1}^4
\hspace{-0.1cm} z_{a} \, (e_{k_a} - e_{-{k_a}})~, \label{unectg} \feq
$$ k_1 := (1,1,1),~~ k_2 := (1,1,-1),~~ k_3 := (1,-1,1),~~ k_4 := (-1,1,1)~; $$
$$ z_{1} := z_2 := (-1,1,0),~ z_{3} := (-1,-1,0),~ z_{4} := -z_3~. $$
The third component of $\uz$ vanishes; however,
this component does not vanish in the exact NS solution $u$ with
this datum
({\footnote{In fact, denoting with $^{(3)}$ the third component we have
$(d u^{(3)}/ d t)(0) = R \PPP(\uz, \uz)^{(3)}$,  which is
nonzero if $R \neq 0$}}), nor in the coefficients $u_1, u_2,...$
of the Reynolds expansion. \par
The above datum was considered by Taylor and Green in
\cite{Tay} for a
pioneering computation of the dissipation rate
of the kinetic energy via a Taylor expansion in time of
the NS solution $u$. The same datum has
been the subject of many subsequent investigations;
among them we cite, in particular, \cite{Bra1}.
These investigations treated
sophisticated issues, such as the numerical verification of Kolmogorov's
hypothesis on turbulence for very large $R$; as already
mentioned, global existence was essential assumed without proof
for these large values of $R$. \par
Eqs.\,\rref{vs} \rref{elles} \rref{defrey} for this datum give
$V_{*} = 1/2$, $L_{*} = 2 \pi/\sqrt{3}$ and
\beq \Rey = {\pi \over \sqrt{3}} \, R = 1.813... \, R~. \feq
The TG symmetries are described in Appendix A; in particular, $\Isott^{+}(\uz)$
has $16$ elements and coincides with $\Isott^{-}(\uz)$. As anticipated, we have used symmetry considerations
to perform the Reynolds expansion up to the order $N=20$. \par
Again for an appreciation of the computational
complexity, we mention that $u_{20}$ has $10560$ nonzero Fourier
coefficients, whose wave vectors are partitioned
in $715 $ orbits under the action of $\Isott^{+}(\uz)$ on $\Zt$.
As an example consider
$u^{(1)}_{20, k}(t)$ for $k = (1,1,1) $, where $^{(1)}$ denotes the first of the three
components; this is a polynomial of degrees $9$ in $t$ and
$547$ in $e^{-t}$.
\par
The expansion up to $N=20$ and its a posteriori analysis give a picture as in
items (i)-(iii) before Eq.\rref{unec0}; for $N=20$ one has \beq R_{\crit} \in
(\rtg, \rtgp), \quad  \mbox{whence} \quad \Rey_{\crit} \in (\reytg, \reytgp)~.
\feq The forthcoming Boxes 3a-3d and 4a-4d present some results of these
computations.
\begin{figure}
\framebox{
\parbox{2in}{
\includegraphics[
height=1.3in,
width=2.0in
]%
{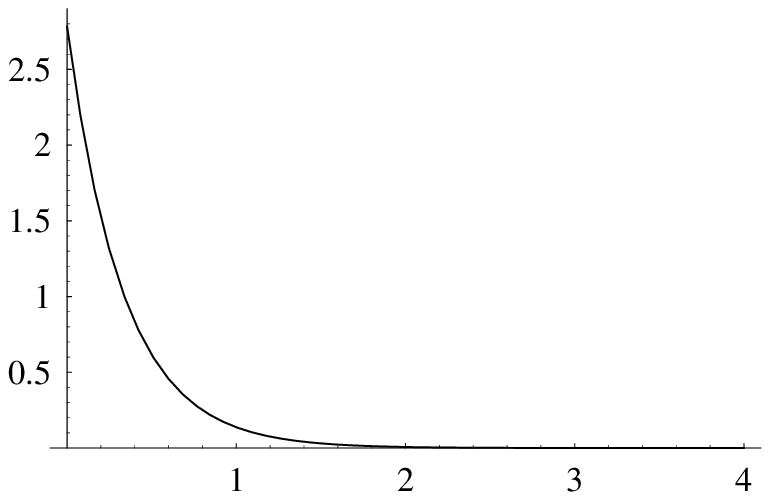}%
\par
{\tiny{
{\textbf{Box 3a.} TG, $N=20$, $R = \rtg$:
the function $\gamma(t) := (2 \pi)^{3/2} |u^{20}_{(1,1,1)}(t)|$. One has
$\gamma(0) = 2.784...$, $\gamma(0.5) = 0.6158...$,
$\gamma(1) = 0.1372...$, $\gamma(1.5) = 0.03061... $, $\gamma(2) = 6.831... \times 10^{-3}$,
$\gamma(4) = 1.693... \times 10^{-5}$,
$\gamma(8) = 1.040... \times 10^{-10}$, $\gamma(10) = 2.579... \times 10^{-13}$~.
}}
\par}
\label{f3a}
}
}
\hskip 0.4cm
\framebox{
\parbox{2in}{
\includegraphics[
height=1.3in,
width=2.0in
]%
{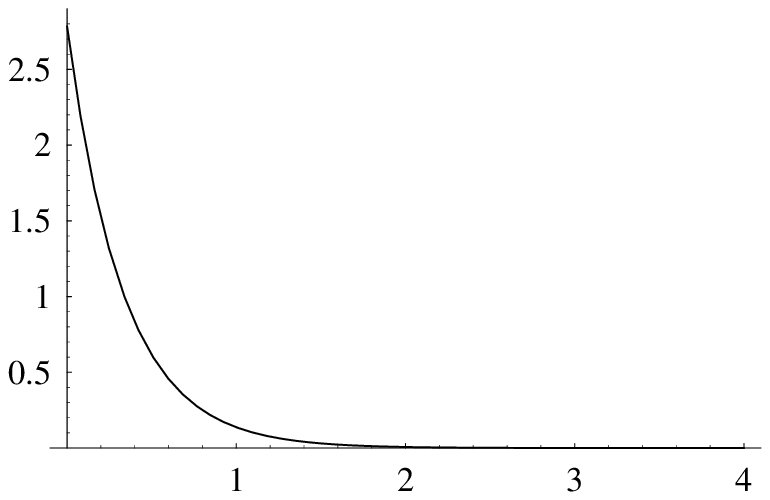}%
\par
{\tiny{
{\textbf{Box 4a.} TG, $N=20$, $R=\rtgp$: the function
$\gamma(t) := (2 \pi)^{3/2} |u^{20}_{(1,1,1)}(t)|$. One has
$\gamma(0) = 2.784...$, $\gamma(0.5) = 0.6154...$,
$\gamma(1) = 0.1371...$, $\gamma(1.5) = 0.03059... $, $\gamma(2) = 6.826... \times 10^{-3}$,
$\gamma(4) = 1.692... \times 10^{-5}$~.
\vskip 0.2cm
{~}
}}
\par}
\label{f4a}
}
}
\framebox{
\parbox{2in}{
\includegraphics[
height=1.3in,
width=2.0in
]%
{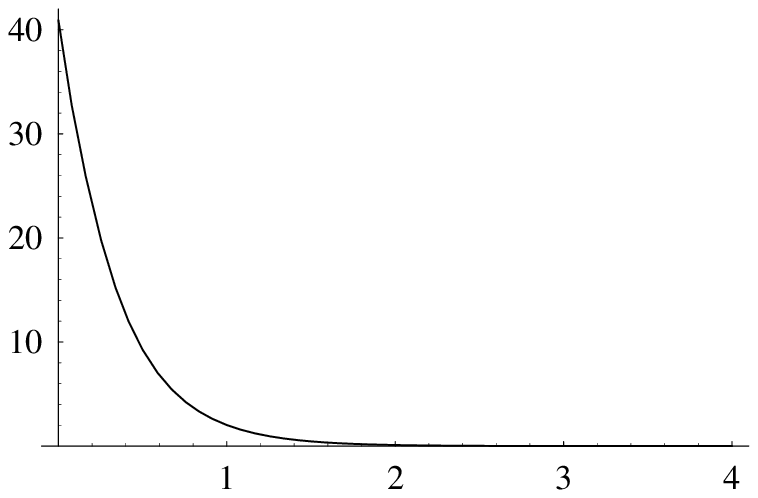}%
\par
{\tiny{
{\textbf{Box 3b.} TG, $N=20$, $R=\rtg$: the function
$\Dd_3(t)$.
One has $\Dd_3(0) = 40.91...$, $\Dd_3(0.5) = 9.257...$,
$\Dd_3(1) = 2.021...$, $\Dd_3(1.5) = 0.4500... $, $\Dd_3(2) = 0.1004...$,
$\Dd_3(4) = 2.488 ... \times 10^{-4}$, $\Dd_3(8) =1.529... \times 10^{-9}$,
$\Dd_3(10) = 3.790... \times 10^{-12}$~.
}}
\par}
\label{f3b}
}
}
\hskip 0.4cm
\framebox{
\parbox{2in}{
\includegraphics[
height=1.3in,
width=2.0in
]%
{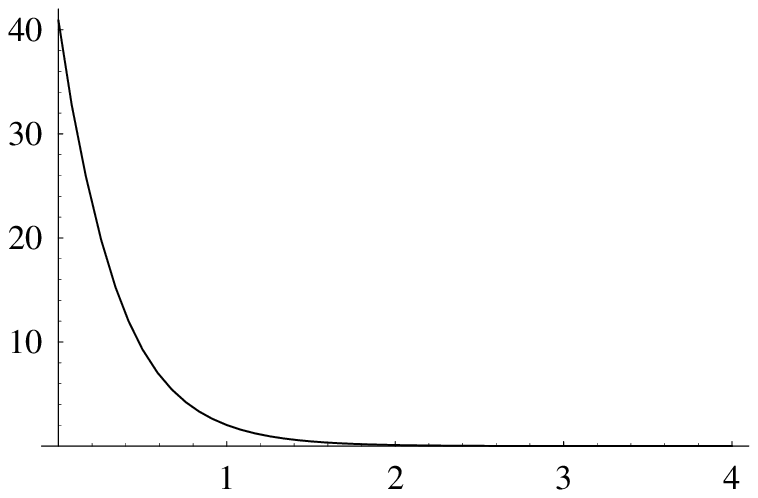}%
\par
{\tiny{
{\textbf{Box 4b.} TG, $N=20$, $R=\rtgp$: the function
$\Dd_3(t)$. One has $\Dd_3(0) = 40.91...$, $\Dd_3(0.5) = 9.266...$,
$\Dd_3(1) = 2.020...$, $\Dd_3(1.5) = 0.4497... $, $\Dd_3(2) = 0.1003...$,
$\Dd_3(4) = 2.487... \times 10^{-4}$, $\Dd_3(8) = 1.528... \times 10^{-9}$,
$\Dd_3(10) = 3.787... \times 10^{-12}$~.
\vskip 0.15cm
{~}
}}
\par}
\label{f4b}
}
}
\framebox{
\parbox{2in}{
\includegraphics[
height=1.3in,
width=2.0in
]%
{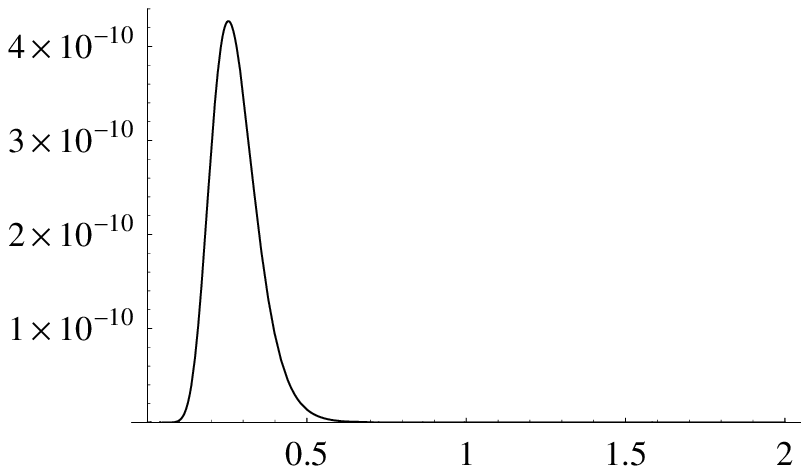}%
\par
{\tiny{
{\textbf{Box 3c.} TG, $N=20$, $R=\rtg$: the function
$\ep_3(t)$.
One has $\ep_3(0) = 0$, $\ep_3(0.25) = 4.263 ... \times 10^{-10}$,
$\ep(0.4) = 9.266...\times 10^{-11}$, $\epsilon(0.6) = 1.661...\times 10^{-12}$,
$\ep_3(1) = 9.152 ... \times 10^{-16}$,
$\ep_3(2) = 6.705... \times 10^{-19}$\, .
}}
\par}
\label{f3c}
}
}
\hskip 0.4cm
\framebox{
\parbox{2in}{
\includegraphics[
height=1.3in,
width=2.0in
]%
{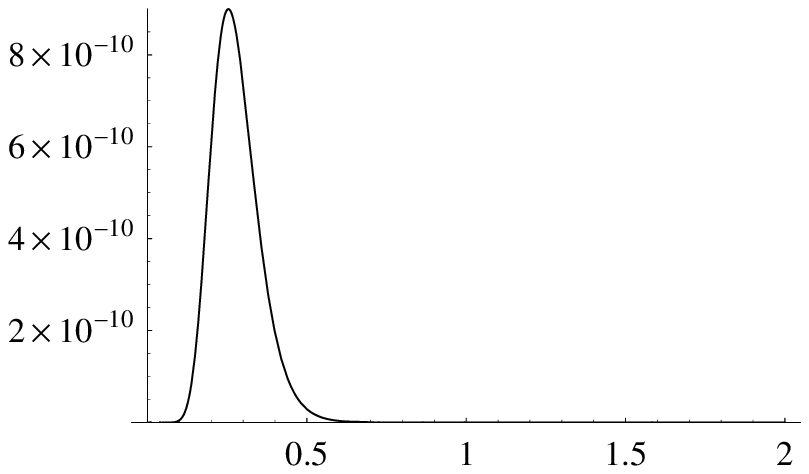}%
\par
{\tiny{
{\textbf{Box 4c.} TG, $N=20$, $R=\rtgp$: the function
$\ep_3(t)$. One has $\ep_3(0) = 0$,
$\ep_3(0.25) = 8.982... \times 10^{-10}$,
$\ep(0.4) = 1.952...\times 10^{-10}$, $\epsilon(0.6) = 3.496...\times 10^{-12}$,
$\ep_3(1) = 1.922...\times 10^{-15}$,
$\ep_3(2) = 1.403... \times 10^{-18}$.}}
\par}
\label{f4c}
}
}
\framebox{
\parbox{2in}{
\includegraphics[
height=1.3in,
width=2.0in
]%
{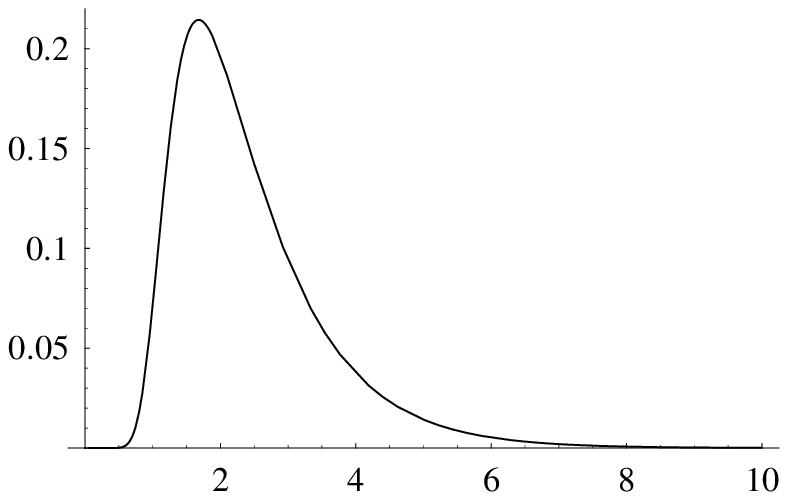}%
\par
{\tiny{
{\textbf{Box 3d.} TG, $N=20$, $R=\rtg$: the function
$\Rr_3(t)$. This
appears to be globally defined,
and vanishing at $+\infty$. One has $\Rr_3(0) = 0$, $\Rr_3(1) = 0.07176...$,
$\Rr_3(1.7) = 0.2143...$, $\Rr_3(2) = 0.1964...$, $\Rr_3(4) = 0.03753...$,
$\Rr_3(8) = 7.202...\times 10^{-4}$, $\Rr_3(10) = 9.754... \times 10^{-5}$~.
\vskip 0.1cm
{~}
}}
\par}
\label{f3d}
}
}
\hskip 3.7cm
\framebox{
\parbox{2in}{
\includegraphics[
height=1.3in,
width=2.0in
]%
{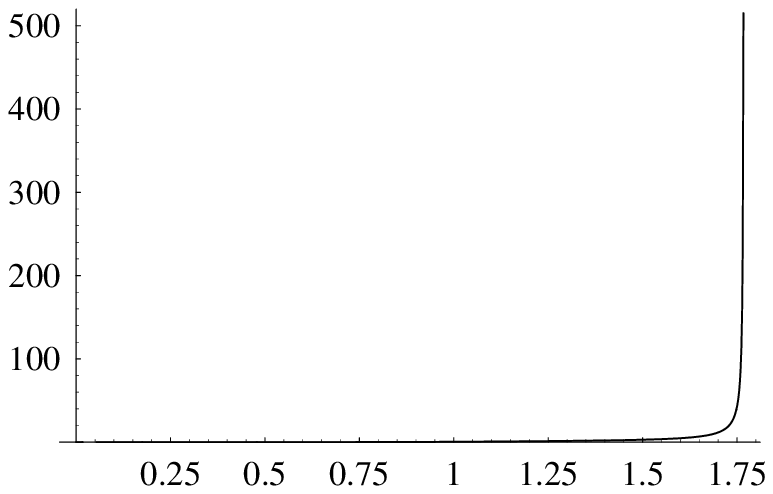}%
\par
{\tiny{
{\textbf{Box 4d.} TG, $N=20$, $R=\rtgp$: the function $\Rr_3(t)$.
This diverges for $t \vain \Tc = 1.768...$\,.
One has $\Rr_3(0) = 0$, $\Rr_3(0.5) = 6.618... \times 10^{-4}$,
$\Rr_3(1) = 0.4435...$, $\Rr_3(1.5) = 2.926...$,
$\Rr_3(1.7) = 11.61...$, $\Rr_3(1.765) = 223.2$~.
\vskip 0.25cm
{~}
}}
\par}
\label{f4d}
}
}
\end{figure}
\vfill\eject\noindent
\textbf{The KM datum.} This is \beq \uz(x_1, x_2, x_3) :=
2 \big( \sin x_1 \cos x_2  \cos x_3 (\cos 2 x_2 - \cos 2 x_3 ) , \label{km}
\feq
$$ \cos x_1  \sin x_2 \cos x_3  (\cos 2 x_3  - \cos 2 x_1 )  ,
\cos x_1  \cos x_2 \sin x_3 (\cos 2 x_1  - \cos 2 x_2 )
\big)~.
$$
Equivalently,
\beq \uz = {i \over 8} \sum_{a=1}^{12} z_a (e_{k_a} - e_{-k_a})~,
\feq
$$ k_1 := (3, 1, 1), \quad k_2 := (3, 1, -1), \quad k_3 := (1, 3, 1), \quad k_4 := (1, 3, -1), $$
$$ k_5 := (1, 1, 3), \quad k_6 := (1, 1, -3)~, \quad k_7 := (1, -1, 3), \quad k_8 := (1, -1, -3), $$
$$ k_9 := (1, -3, 1), \quad k_{10} := (1, -3, -1), \quad k_{11} := (3, -1, 1), \quad k_{12} := (3, -1, -1), $$
$$ z_1 := (0,1,-1), \quad z_2 := (0,1,1), \quad z_3 := z_9 := (-1,0,1), \quad z_4 := z_{10} := (-1,0,-1)~, $$
$$ z_5 := z_6 := (1,-1,0), \quad z_7 := z_8 := (1,1,0)~, \quad z_{11} := - z_2, \quad z_{12} := -z_1~. $$
This datum was the subject of an investigation started
by Kida \cite{Kid} and continued by Kida and Murakami \cite{Kid2};
it is maximally symmetric, in the sense that
$\Isott^{+}(\uz)$ is the full octahedral group $O(3,\interi)$,
and coincides with $\Isott^{-}(\uz)$  (see Appendix A).
In the cited works, this feature was used to reduce the computational costs
in the solution of the NS equations via pseudo-spectral methods
(with no discussion of the global existence problem,
as typical of numerical investigations on turbulence).
\par
Eqs.\,\rref{vs} \rref{elles} \rref{defrey} for this datum give
$V_{*} = \sqrt{3}/2$,  $L_{*} = 2 \pi/\sqrt{11}$ and
\beq \Rey = \sqrt{{3 \over 11}}~ \pi \, R = 1.640... \, R~. \feq
In the KM case, using the symmetries
we could perform the Reynolds expansion up to the order $N=12$.
Let us mention that $u_{12}$ has $33312 $ nonzero Fourier coefficients,
whose wave vectors are partitioned
in $797$ orbits under the action of $\Isott^{+}(\uz) = O(3, \interi)$ on $\Zt$.
As an example consider
$u^{(2)}_{12, k}(t)$ for $k = (3,1,1) $, where $^{(2)}$ denotes the second of the three
components; this is a polynomial of degrees $5$ in $t$ and
$867$ in $e^{-t}$.
\par
The result of computations up to $N=12$ is a picture as in items (i)-(iii)
before \rref{unec0}; for $N=12$ one has \beq R_{\crit} \in (\rkm, \rkmp), \quad
\mbox{whence} \quad \Rey_{\crit} \in (\reykm, \reykmp)~. \feq
The forthcoming
Boxes 5a-5d and 6a-6d present some results of these computations. (Note
that the numerical sample values reported in Boxes 5a and 6a, 5b and 6b
are almost always the same to 4 meaningful digits; the situation is
different for Boxes 5c and 6c, 5d and 6d.)
\begin{figure}
\framebox{
\parbox{2in}{
\includegraphics[
height=1.3in,
width=2.0in
]%
{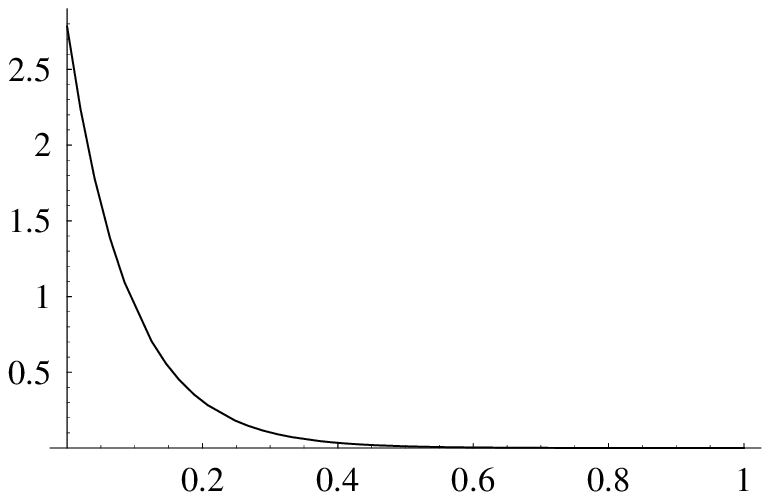}%
\par
{\tiny{
{\textbf{Box 5a.} KM, $N=12$, $R = \rkm$: the function
$\gamma(t) := (2 \pi)^{3/2} |u^{12}_{(3,1,1)}(t)|$. One has
$\gamma(0) = 2.784...$, $\gamma(0.5) = 0.01137...$,
$\gamma(1) = 4.648...\times 10^{-5}$, $\gamma(1.5) = 1.899... \times 10^{-7}$,
$\gamma(2) = 7.762...\times 10^{-10}$,
$\gamma(3) = 1.296...\times 10^{-14}$,
$\gamma(4) = 2.165... \times 10^{-19}$~.
}}
\par}
\label{f5a}
}
}
\hskip 0.4cm
\framebox{
\parbox{2in}{
\includegraphics[
height=1.3in,
width=2.0in
]%
{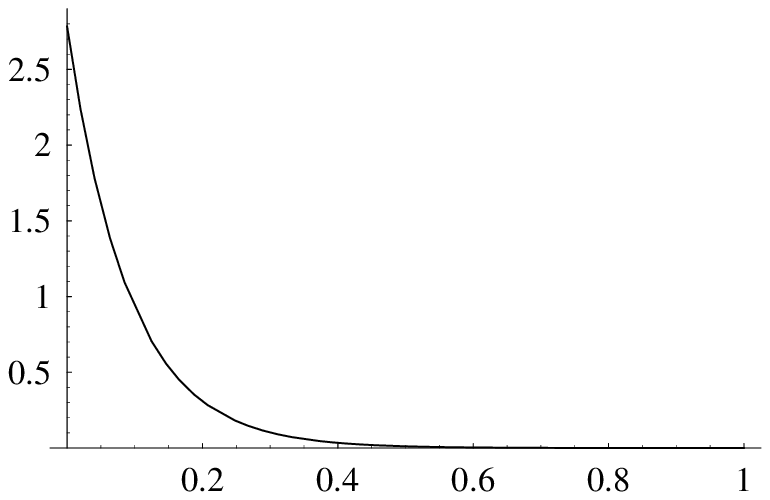}%
\par
{\tiny{
{\textbf{Box 6a.} KM, $N=12$, $R=\rkmp$: the function
$\gamma(t) := (2 \pi)^{3/2} |u^{12}_{(3,1,1)}(t)|$. One has
$\gamma(0) = 2.784...$, $\gamma(0.5) = 0.01137...$,
$\gamma(1) = 4.647...\times 10^{-5}$, $\gamma(1.5) = 1.899... \times 10^{-7}$,
$\gamma(2) = 7.762...\times 10^{-10}$,
$\gamma(3) = 1.296...\times 10^{-14}$,
$\gamma(4) = 2.165... \times 10^{-19}$~.
\vskip 0.15cm
{~}
}}
\par}
\label{f6a}
}
}
\framebox{
\parbox{2in}{
\includegraphics[
height=1.3in,
width=2.0in
]%
{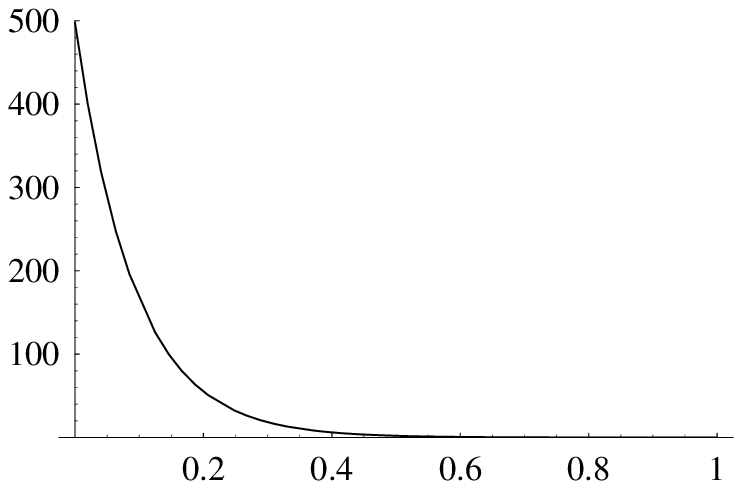}%
\par
{\tiny{
{\textbf{Box 5b.} KM, $N=12$, $R=\rkm$: the function
$\Dd_3(t) = \| u^{12}(t) \|_3$.
One has $\Dd_3(0) = 497.6...$, $\Dd_3(0.5) = 2.032...$,
$\Dd_3(1) = 8.307 ...\times 10^{-3}$,
$\Dd_3(1.5) = 3.395... \times 10^{-5}$, $\Dd_3(2) = 1.387...\times 10^{-7}$,
$\Dd_3(3) = 2.317...\times 10^{-12}$,
$\Dd_3(4) = 3.870...\times 10^{-17}$~.
}}
\par}
\label{f5b}
}
}
\hskip 0.4cm
\framebox{
\parbox{2in}{
\includegraphics[
height=1.3in,
width=2.0in
]%
{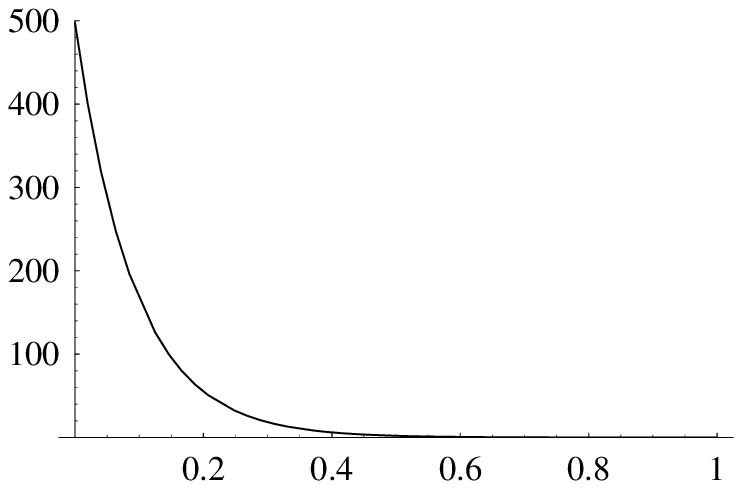}%
\par
{\tiny{
{\textbf{Box 6b.} KM, $N=12$, $R=\rkmp$: the function
$\Dd_3(t) = \| u^{12}(t) \|_3$.
One has $\Dd_3(0) = 497.6...$, $\Dd_3(0.5) = 2.032...$,
$\Dd_3(1) = 8.307 ...\times 10^{-3}$,
$\Dd_3(1.5) = 3.394... \times 10^{-5}$, $\Dd_3(2) = 1.387...\times 10^{-7}$,
$\Dd_3(3) = 2.317...\times 10^{-12}$,
$\Dd_3(4) = 3.870...\times 10^{-17}$~.
\vskip 0.15cm
{~}
}}
\par}
\label{f6b}
}
}
\framebox{
\parbox{2in}{
\includegraphics[
height=1.3in,
width=2.0in
]%
{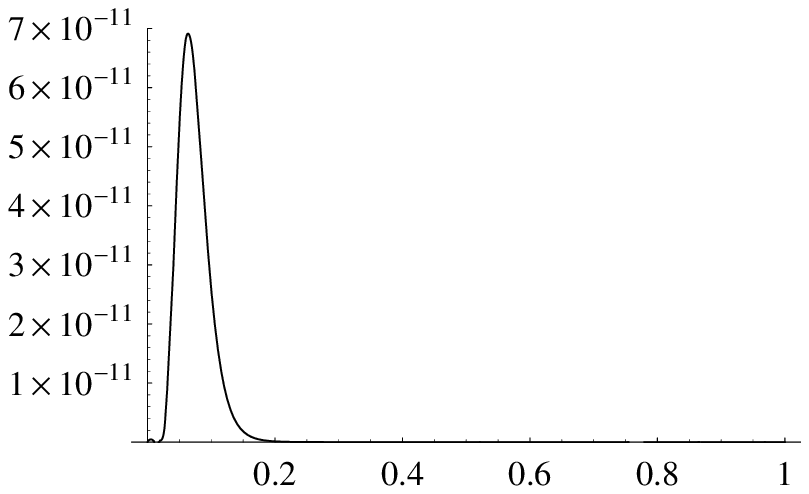}%
\par
{\tiny{
{\textbf{Box 5c.} KM, $N=12$, $R=\rkm$: the function
$\ep_3(t)$. One has $\ep_3(0) = 0$,
$\ep_3(0.06) = 6.820... \times 10^{-11}$,
$\ep_3(0.1) = 2.577... \times 10^{-11}$,
$\ep_3(0.5) = 3.420... \times 10^{-18}$,
$\ep_3(2) = 6.053... \times 10^{-33}$,
$\ep_3(4) = 4.710... \times 10^{-52}$\, .
}}
\par}
\label{f5c}
}
}
\hskip 0.4cm
\framebox{
\parbox{2in}{
\includegraphics[
height=1.3in,
width=2.0in
]%
{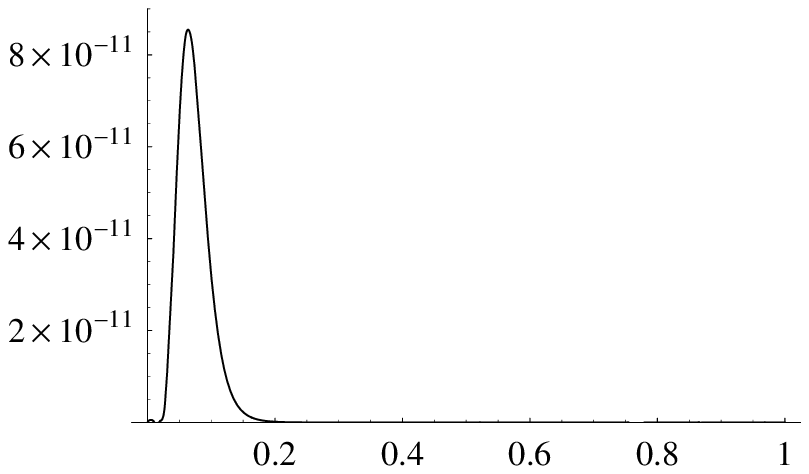}%
\par
{\tiny{
{\textbf{Box 6c.} KM, $N=12$, $R=\rkmp$: the function
$\ep_3(t)$. One has $\ep_3(0) = 0$,
$\ep_3(0.06) = 8.432... \times 10^{-11}$,
$\ep_3(0.1) = 3.187... \times 10^{-11}$,
$\ep_3(0.5) = 4.226... \times 10^{-18}$,
$\ep_3(2) = 7.478... \times 10^{-33}$,
$\ep_3(4) = 5.819... \times 10^{-52}$\, .
\vskip 0.05cm
{~}
}}
\par}
\label{f6c}
}
}
\framebox{
\parbox{2in}{
\includegraphics[
height=1.3in,
width=2.0in
]%
{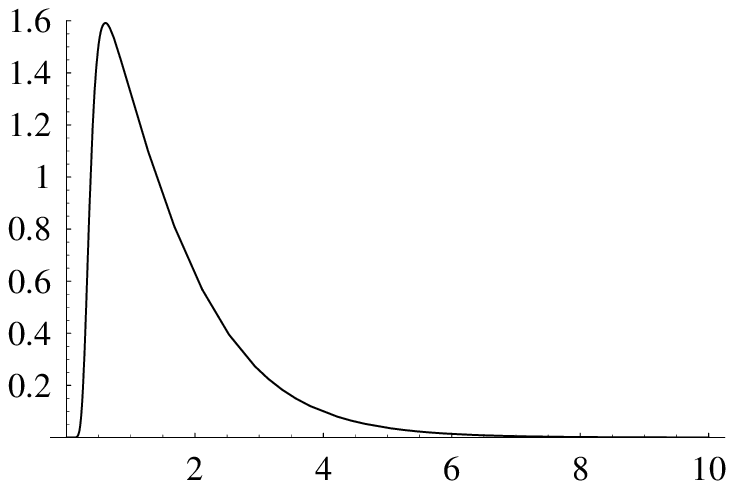}%
\par
{\tiny{
{\textbf{Box 5d.} KM, $N=12$, $R=\rkm$: the function $\Rr_3(t)$. This
appears to be \hbox{globally} defined, and vanishing at $+\infty$.
One has $\Rr_3(0) = 0$, $\Rr_3(0.3) = 0.4585...$,
$\Rr_3(0.6) = 1.591...$,
$\Rr_3(1) = 1.319....$, $\Rr_3(2) = 0.6250...$, $\Rr_3(4) = 0.09886...$,
$\Rr_3(8) = 1.859...\times 10^{-3}$, $\Rr_3(10) = 2.517... \times 10^{-4}$~.
}}
\par}
\label{f5d}
}
}
\hskip 3.7cm
\framebox{
\parbox{2in}{
\includegraphics[
height=1.3in,
width=2.0in
]%
{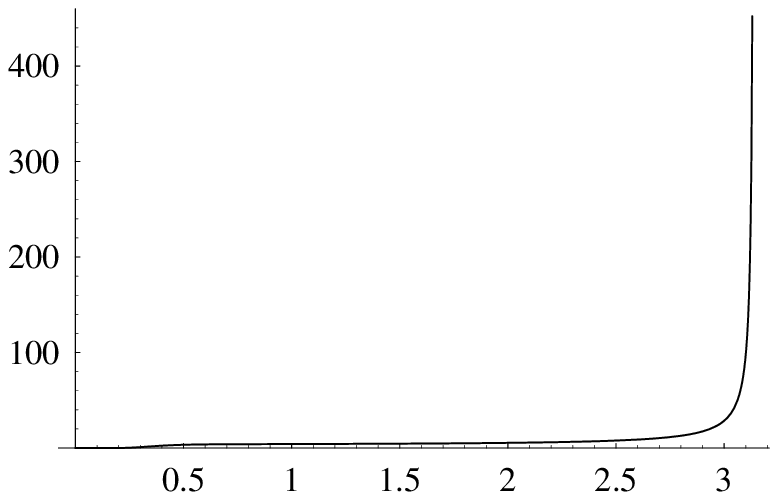}%
\par
{\tiny{
{\textbf{Box 6d.} KM, $N=12$, $R=\rkmp$: the function $\Rr_3(t)$. This diverges
for $t \vain \Tc = 3.138...$\,.
One has $\Rr_3(0) = 0$, $\Rr_3(0.1) = 5.419...\times 10^{-6}$,
$\Rr_3(0.2) =  0.0485...$, $\Rr_3(0.5) =  3.399...$,
$\Rr_3(1) = 4.171...$, $\Rr_3(2) = 5.418...$,
$\Rr_3(3) = 28.53...$, $\Rr_3(3.13) = 452.2...$~.
\vskip 0.25cm
{~}
}}
\par}
\label{f6d}
}
}
\end{figure}
\vfill\eject\noindent
\section{Concluding remarks}
\label{qua}
Let us propose a question that, in our opinion, is worthy
of future consideration: keeping the general setting of the present
paper, is it possible to improve the critical Reynolds
numbers yielding global existence for the initial data considered here?  \par
An obvious attempt one could make in this direction
is to try higher order Reynolds expansions
by means of more powerful computational utilities;
however, it is not granted that this strategy would
yield significant improvements ({\footnote{It might happen
that, for large $R$ or large times, there
is an optimal order $N$ giving the best approximation
of the exact NS solution and that,
for larger $N$, the norms of $u^N(t)$ and of
its differential error increase, finally yielding worse results in the
application of the control Cauchy problem. One could
expect this to happen for values of $R$ or times $t$
for which the series $\sum_{j=0}^{+\infty} R^j u_j(t)$
is not convergent; as already mentioned in \cite{apprey}, the
convergence of this series is known for small $R$ or short
times, but bounds on $R$ or $t$ yielding convergence
are not presently known with sufficient precision.}}). \par
The problem of global existence for
higher Reynolds numbers could be attacked
via a different strategy.
In this case the idea is to devise
specific versions of the control Cauchy problem \rref{cont}, fitted
to the symmetries of the initial data under investigation.
In particular, one could consider the basic inequality \rref{basineq} and the Kato inequality
\rref{katineq} in the subspaces of $\HM{n}$ and $\HM{n+1}$
formed by the vector fields which have
the same symmetry group as the initial
datum; the constants $K_n, G_n$ for the inequalities \rref{basineq}
\rref{katineq} in these subspaces could be significantly smaller,
and the control problem \rref{cont} with these smaller constants
would yield global existence for higher Reynolds numbers.\par
However, estimating the constants for
\rref{basineq} \rref{katineq}
in the presence of symmetries requires some
effort: one must adapt the approach
of \cite{cog} \cite{cok} to
the case where a symmetry group
is specified and, especially,
one must make anew rather expensive
numerical computations to evaluate
$K_n$ and $G_n$ for specific values
of $n$ and for the given symmetries.
We plan to treat the above issues in future works.
\vskip 0.6cm \noindent
\textbf{Acknowledgments.}
This work was partly supported by INdAM, INFN and by MIUR, PRIN 2010
Research Project  ``Geometric and analytic theory of Hamiltonian systems in finite and infinite dimensions''.
\vfill \eject \noindent
\appendix
\section{Appendix. Symmetries of the BNW, TG and KM data}
Finding the symmetries or pseudo-symmetries of any NS initial datum
$\uz$ amounts to determine all pairs $(S,a) \in \ohz$ such that
$\E_{*}(S,a) \uz = \pm \uz$. When $\uz$ is a Fourier polynomial it is
generally convenient to rephrase this equation
in terms of Fourier coefficients via \rref{pushf};
the solutions can be obtained by automatic computations,
say with Mathematica.
These remarks apply, in particular, to the three initial
data considered in this paper.
\par
Throughout this Appendix we work in dimension $d=3$, using the following notations:
$$
\barray{llll}
D_1 := \mbox{diag}(1,1,1), &
D_2 := \mbox{diag}(-1,1,1), &
D_3 := \mbox{diag}(1,-1,1), &
D_4 := \mbox{diag}(1,1,-1), \\
D_5 := - D_2, &
D_6 := - D_3, &

D_7 := - D_4, &
D_8 := - D_1
\farray
$$
$$
Q_1 := \left( \barray{ccc} 1 & 0 & 0 \\  0 & 1 & 0 \\ 0 & 0 & 1 \farray \right) ,~~
Q_2 := \left( \barray{ccc} 0 & 1 & 0 \\  0 & 0 & 1 \\ 1 & 0 & 0 \farray \right) ,~~
Q_3 := \left( \barray{ccc} 0 & 0 & 1 \\  1 & 0 & 0 \\ 0 & 1 & 0 \farray \right) ,
$$
\beq
{~} \hspace{0.7cm}
Q_4 := \left( \barray{ccc} 0 & 1 & 0 \\  1 & 0 & 0 \\ 0 & 0 & 1 \farray \right) ,~~
Q_5 := \left( \barray{ccc} 1 & 0 & 0 \\  0 & 0 & 1 \\ 0 & 1 & 0 \farray \right) ,~~
Q_6 := \left( \barray{ccc} 0 & 0 & 1 \\  0 & 1 & 0 \\ 1 & 0 & 0 \farray \right) .
\feq
\normalsize
Any matrix $Q_\beta$ ($\beta=1,...,6$) acts on elements of $\reali^3$ or $\Tt$ applying
to their components one of the $6$ permutations of $\{1,2,3\}$. According
to Eq.\,\rref{rep}, the octahedral group $O(3, \interi)$ is formed by all matrices
of the form
\beq S_{\alpha \beta} := D_{\alpha} Q_{\beta} \qquad (\alpha=1,...,8;\,  \beta = 1,...,6)~. \feq
To go on, we put
\beq
a_1 := (0,0,0),~~
a_2 := (\pi,0,0),~~
a_3 := (0,\pi,0),~~
a_4 := (0,0,\pi), \feq
$$ a_5 := (\pi,\pi,0),~~
a_6 := (\pi,0,\pi),~~
a_7 := (0,\pi,\pi),~~
a_8 := (\pi,\pi,\pi)
$$
(where $\pi$ is an abbreviation for $\pi$ mod.$\,2 \pi \interi$).
With the above notations, the symmetries and pseudosymmetries of the BNW, TG and KM data
can be described as follows.
\vfill \eject \noindent
\textbf{BNW case.} Let $\uz$ denote the BNW datum \rref{unec0}. Then \parn
\vbox{
\beq \Isot^{+}(\uz) = \{ (S_{\alpha \beta}, a_{\gamma})~|~(\alpha, \beta, \gamma) \in I \}~, \feq
\footnotesize
$$
{~}
\hskip 0.6cm
\barray{llllllll}
\! \! \! \! \! \! \! \! \! \! \! \! \! \! \! \! \! \!
I :=
\{~ (1, 1, 1), &  (1, 1, 8), &  (1, 2, 4), & (1, 2, 5), & (1, 3, 2), & (1, 3, 7), \\
(8, 4, 4), & (8, 4, 5), & (8, 5, 2), & (8, 5, 7), & (8, 6, 1), & (8, 6, 8)~ \}~.
\farray
$$
}
\normalsize
This group has 12 elements; it was already described (with different
notations) in \cite{bnw}, where it was shown to be isomorphic
to the dihedral group $\textbf{D}_6$ (the group of orthogonal symmetries
of a regular hexagon). The pseudo-symmetry space of the BNW datum
contains $(S_{8 1}, a_1) = (\mbox{diag}(-1,-1,-1), (0,0,0))$;
thus
\beq \Isot^{-}(\uz) = \{(S_{\alpha \beta}, a_{\gamma}) \circ (S_{8 1}, a_1)~|~(\alpha, \beta, \gamma) \in I \}
= \{ (-S_{\alpha \beta}, a_{\gamma})~|~(\alpha, \beta, \gamma) \in I \}~. \feq
The reduced symmetry group
and pseudo-symmetry space for this datum are
\beq \Isott^{\pm}(\uz) = \{ \pm S_{\alpha \beta}~|~(\alpha, \beta, \gamma) \in I \}~; \feq
the above two sets are disjoint, and each one of them has $6$ elements.
In \cite{bnw}, it was shown that $\Isott^{+}(\uz)$ is isomorphic
to the dihedral group $\textbf{D}_3$ (the group of orthogonal symmetries
of an equilateral triangle).
\vskip 0.2cm \noindent
\textbf{TG case.} Let $\uz$ denote the TG datum \rref{tg}. Then \parn
\vbox{
\beq \Isot^{+}(\uz) = \{ (S_{\alpha \beta}, a_{\gamma})~|~(\alpha, \beta, \gamma) \in I \}~, \feq
\footnotesize
$$
{~}
\hskip 0.6cm
\barray{llllllll}
\! \! \! \! \! \! \! \! \! \! \! \! \! \! \! \! \! \!
I := \{~(1, 1, 1), &  (1, 1, 5), &  (1, 1, 6), &  (1, 1, 7), &  (1, 4, 2), &  (1, 4, 3), &  (1, 4, 4), &  (1, 4, 8), \\
(2, 1, 1), &  (2, 1, 5), &  (2, 1, 6), &  (2, 1, 7), &  (2, 4, 2), &  (2, 4, 3), &  (2, 4, 4), &  (2, 4, 8), \\
 (3, 1, 1), &  (3, 1, 5), &  (3, 1, 6), &  (3, 1, 7), &  (3, 4, 2), &  (3, 4, 3), &  (3, 4, 4), &  (3, 4, 8), \\
 (4, 1, 1), &  (4, 1, 5), &  (4, 1, 6), &  (4, 1, 7), &   (4, 4, 2), &  (4, 4, 3), &  (4, 4, 4), &  (4, 4, 8), \\
 (5, 1, 1), &  (5, 1, 5), &  (5, 1, 6), &  (5, 1, 7), &  (5, 4, 2), &  (5, 4, 3), &  (5, 4, 4), &  (5, 4, 8), \\
 (6, 1, 1), &  (6, 1, 5), &  (6, 1, 6), &  (6, 1, 7), &  (6, 4, 2), &  (6, 4, 3), &  (6, 4, 4), &  (6, 4, 8), \\
 (7, 1, 1), &  (7, 1, 5), &  (7, 1, 6), &  (7, 1, 7), &  (7, 4, 2), &  (7, 4, 3), &  (7, 4, 4), &  (7, 4, 8), \\
 (8, 1, 1), &  (8, 1, 5), &  (8, 1, 6), &  (8, 1, 7), &  (8, 4, 2), &  (8, 4, 3), &  (8, 4, 4), &  (8, 4, 8)~\}.
\farray
$$
}
\normalsize
This group has 64 elements. The pseudo-symmetry space of the TG datum
contains $(S_{1 1}, a_8) = (\mbox{diag}(1,1,1), (\pi,\pi,\pi))$; thus
\beq {~} \hspace{-0.4cm} \Isot^{-}(\uz) = \{ (S_{1 1}, a_8) \circ
(S_{\alpha \beta}, a_{\gamma})~|~(\alpha, \beta, \gamma) \in I \} =
\{(S_{\alpha \beta}, a_8 + a_{\gamma})~|~(\alpha, \beta, \gamma) \in I \}.
\feq
The reduced symmetry group
and pseudo-symmetry space for this datum coincide; they have
16 elements, and are given by
\beq \Isott^{\pm}(\uz) = \{ S_{\alpha \beta}~|~(\alpha, \beta, \gamma) \in I \}~. \feq
\vfill \eject \noindent
\textbf{KM case.} Let $\uz$ denote the KM datum \rref{km}. Then \parn
\vbox{
\beq \Isot^{+}(\uz) = \{ (S_{\alpha \beta}, a_{\gamma})~|~(\alpha, \beta, \gamma) \in I \}~, \feq
\footnotesize
$$
{~}
\hskip 0.6cm
\barray{llllllll}
\! \! \! \! \! \! \! \! \! \! \! \! \! \! \! \! \! \!
I := \{~ (1, 1, 1), & (1, 1, 5), & (1, 1, 6), & (1, 1, 7), & (1, 2, 1), & (1, 2, 5), & (1, 2, 6), & (1, 2, 7), \\
(1, 3, 1), & (1, 3, 5), & (1, 3, 6), & (1, 3, 7), & (1, 4, 2), & (1, 4, 3), & (1, 4, 4), & (1, 4, 8), \\
(1, 5, 2), & (1, 5, 3), & (1, 5, 4), & (1, 5, 8), & (1, 6, 2), & (1, 6, 3), & (1, 6, 4), & (1, 6, 8), \\
(2, 1, 1), & (2, 1, 5), & (2, 1, 6), & (2, 1, 7), & (2, 2, 1), & (2, 2, 5), & (2, 2, 6), & (2, 2, 7), \\
(2, 3, 1), & (2, 3, 5), & (2, 3, 6), & (2, 3, 7), & (2, 4, 2), & (2, 4, 3), & (2, 4, 4), & (2, 4, 8), \\
(2, 5, 2), & (2, 5, 3), & (2, 5, 4), & (2, 5, 8), & (2, 6, 2), & (2, 6, 3), & (2, 6, 4), & (2, 6, 8), \\
(3, 1, 1), & (3, 1, 5), & (3, 1, 6), & (3, 1, 7), & (3, 2, 1), & (3, 2, 5), & (3, 2, 6), & (3, 2, 7), \\
(3, 3, 1), & (3, 3, 5), & (3, 3, 6), & (3, 3, 7), & (3, 4, 2), & (3, 4, 3), & (3, 4, 4), & (3, 4, 8), \\
(3, 5, 2), & (3, 5, 3), & (3, 5, 4), & (3, 5, 8), & (3, 6, 2), & (3, 6, 3), & (3, 6, 4), & (3, 6, 8), \\
(4, 1, 1), & (4, 1, 5), & (4, 1, 6), & (4, 1, 7), & (4, 2, 1), & (4, 2, 5), & (4, 2, 6), & (4, 2, 7), \\
(4, 3, 1), & (4, 3, 5), & (4, 3, 6), & (4, 3, 7), & (4, 4, 2), & (4, 4, 3), & (4, 4, 4), & (4, 4, 8), \\
(4, 5, 2), & (4, 5, 3), & (4, 5, 4), & (4, 5, 8), & (4, 6, 2), & (4, 6, 3), & (4, 6, 4), & (4, 6, 8), \\
(5, 1, 1), & (5, 1, 5), & (5, 1, 6), & (5, 1, 7), & (5, 2, 1), & (5, 2, 5), & (5, 2, 6), & (5, 2, 7), \\
(5, 3, 1), & (5, 3, 5), & (5, 3, 6), & (5, 3, 7), & (5, 4, 2), & (5, 4, 3), & (5, 4, 4), & (5, 4, 8), \\
(5, 5, 2), & (5, 5, 3), & (5, 5, 4), & (5, 5, 8), & (5, 6, 2), & (5, 6, 3), & (5, 6, 4), & (5, 6, 8), \\
(6, 1, 1), & (6, 1, 5), & (6, 1, 6), & (6, 1, 7), & (6, 2, 1), & (6, 2, 5), & (6, 2, 6), & (6, 2, 7), \\
(6, 3, 1), & (6, 3, 5), & (6, 3, 6), & (6, 3, 7), & (6, 4, 2), & (6, 4, 3), & (6, 4, 4), & (6, 4, 8), \\
(6, 5, 2), & (6, 5, 3), & (6, 5, 4), & (6, 5, 8), & (6, 6, 2), & (6, 6, 3), & (6, 6, 4), & (6, 6, 8), \\
(7, 1, 1), & (7, 1, 5), & (7, 1, 6), & (7, 1, 7), & (7, 2, 1), & (7, 2, 5), & (7, 2, 6), & (7, 2, 7), \\
(7, 3, 1), & (7, 3, 5), & (7, 3, 6), & (7, 3, 7), & (7, 4, 2), & (7, 4, 3), & (7, 4, 4), & (7, 4, 8), \\
(7, 5, 2), & (7, 5, 3), & (7, 5, 4), & (7, 5, 8), & (7, 6, 2), & (7, 6, 3), & (7, 6, 4), & (7, 6, 8), \\
(8, 1, 1), & (8, 1, 5), & (8, 1, 6), & (8, 1, 7), & (8, 2, 1), & (8, 2, 5), & (8, 2, 6), & (8, 2, 7), \\
(8, 3, 1), & (8, 3, 5), & (8, 3, 6), & (8, 3, 7), & (8, 4, 2), & (8, 4, 3), & (8, 4, 4), & (8, 4, 8), \\
(8, 5, 2), & (8, 5, 3), & (8, 5, 4), & (8, 5, 8), & (8, 6, 2), & (8, 6, 3), & (8, 6, 4), & (8, 6, 8)~ \}.
\farray
$$
}
\normalsize
This group has 192 elements. As in the TG case,
the pseudo-symmetry space of the KM datum
contains $(S_{1 1}, a_8) = (\mbox{diag}(1,1,1), (\pi,\pi,\pi))$; thus
\beq {~} \hspace{-0.3cm} \Isot^{-}(\uz) = \{ (S_{1 1}, a_8) \circ
(S_{\alpha \beta}, a_{\gamma}) | (\alpha, \beta, \gamma) \in I \} =
\{(S_{\alpha \beta}, a_8 + a_{\gamma}) | (\alpha, \beta, \gamma) \in I \}.
\feq
The reduced symmetry group
and pseudo-symmetry space for the KM datum coincide; they have
48 elements, i.e., they coincide with the full octahedral group:
\beq \Isott^{\pm}(\uz) = \{ S_{\alpha \beta}~|~(\alpha, \beta) \in I \} =
O(3,\interi)~. \feq
\vfill \eject \noindent
{~}


\vskip 0cm \noindent
\begin{thebibliography}{99}
\footnotesize{
\bibitem{BKM} J.T. Beale, T. Kato, A.J. Majda, \textsl{Remarks on the breakdown of smooth solutions
for the 3D Euler equations}, Commun. Math. Phys. \textbf{94} (1984), 61-66.
\bibitem{Nec}
E. Behr, J. Ne$\check{\mbox{c}}$as, H. Wu, \textsl{On blow-up of
solution for Euler equations}, M2AN: Math. Model. Numer. Anal.
\textbf{35} (2001), 229-238. \vsm
\bibitem{Bra1} M.E. Brachet, D. Meiron, S. Orszag, B. Nickel,
R. Morf, U. Frisch, \textsl{Small scale structure of the Taylor-Green vortex},
J. Fluid Mech. \textbf{130} (1983), 411-452.
\vsm
\bibitem{CheG} J.Y. Chemin. I. Gallagher,
\textsl{On the global wellposedness of the $3$-D Navier-Stokes
equations with large initial data}, Ann. Sc. Ecole Norm. Sup.
\textbf{39} (2006), 679-698.
\vsm
\bibitem{Che} S.I. Chernyshenko, P. Constantin, J.C. Robinson, E.S. Titi,
\textsl{A posteriori regularity of the three-dimensional
Navier-Stokes equations from numerical computations}, J. Math.
Phys. \textbf{48} (2007), 065204/10. \vsm
\bibitem{Rob3} M. Dashti, J.C. Robinson,
\textsl{An a posteriori condition on the numerical approximations
of the Navier-Stokes equations for the existence of a strong
solution}, SIAM J. Numer. Anal. \textbf{46} (2008), 3136-3150.
\vsm
\bibitem{Gig} Y. Giga,
\textsl{Solutions for semilinear parabolic equations in $L^p$ and regularity of weak solutions of the
Navier-Stokes system}, J. Differential Equations \textbf{62} (1986), 186-212.
\bibitem{Kato} T.Kato, \textsl{Nonstationary flows of viscous and ideal fluids in $\reali^3$},
J.Funct.Anal. \textbf{9} (1972), 296-305. \vsm
\bibitem{Kat2} T. Kato, \textsl{Quasi-linear equations of evolution, with applications to
partial differential equations}, in ``Spectral theory and
differential equations'', Proceedings of the Dundee Symposium,
Lecture Notes in Mathematics \textbf{448} (1975), 23-70. \vsm
\bibitem{Kid}
S. Kida, \textsl{Three-dimensional periodic flows with
high-symmetry}, J. Phys. Soc. Japan \textbf{54} (1985), 2132-2140.
\vsm
\bibitem{Kid2}
S. Kida, Y. Murakami, \textsl{Kolmogorov's
spectrum in a freely decaying turbulence},
J. Phys. Soc. Japan \textbf{55} (1986), 9-12.
\vsm
\bibitem{Koz} H. Kozono, Y. Taniuchi, \textsl{Limiting case of the Sobolev inequality in BMO,
with application to the Euler equations}, Commun. Math. Phys. \textbf{214} (2000), 191--200.
\bibitem{Koz2} H. Kozono, Y. Taniuchi, \textsl{Bilinear estimates in BMO,
and the Navier-Stokes equations}, Math. Z. \textbf{235} (2000), 173--194.
\bibitem{Kuk} I. Kukavica, W. Rusin, M. Ziane,
\textsl{A class of solutions of the Navier-Stokes equations
with large data}, J. Differential Equations \textbf{255} (2013), 1492--1514.
\bibitem{bnw}
C. Morosi, M. Pernici, L. Pizzocchero, \textsl{On power series
solutions for the Euler equation, and the
Behr-Ne$\check{\mbox{c}}$as-Wu initial datum}, ESAIM Math. Model.
Numer. Anal. \textbf{47} (2013), 663-688.
\vsm
\bibitem{padova} C. Morosi, M. Pernici, L. Pizzocchero, \textsl{A posteriori estimates
for Euler and Navier-Stokes equations}, in: F. Ancona, A. Bressan, P. Marcati,
A. Marson (Eds.), Hyperbolic Problems:
Theory, Numerics and Applications, Proceedings
of the XIV International Conference held in Padova  (June 25-29, 2012), in:
AIMS Series on Applied Mathematics \textbf{8} (2014), 847-855.
\vsm
\bibitem{coga} C. Morosi, M. Pernici, L. Pizzocchero,
\textsl{On the constants in some inequalities for the Navier-Stokes quadratic
nonlinearity}, in preparation.
\vsm
\bibitem{uno} C. Morosi, L. Pizzocchero, \textsl{On approximate solutions of semilinear
evolution equations}, Rev. Math. Phys. \textbf{16} (2004),
383-420.
\vsm
\bibitem{due} C. Morosi, L. Pizzocchero, \textsl{On approximate solutions of
semilinear evolution equations II. Generalizations, and
applications to Navier-Stokes equations}, Rev. Math. Phys.
\textbf{20} (2008), 625-706.  \vsm
\bibitem{accau} C. Morosi, L. Pizzocchero,
\textsl{An $H^1$ setting for the Navier-Stokes equations:
Quantitative estimates}, Nonlinear Anal. \textbf{74} (2011),
2398-2414.  \vsm
\bibitem{cog} C. Morosi, L. Pizzocchero, \textsl{On the constants
in a Kato inequality for the Euler and NS equations}, Commun. Pure
Appl. \hbox{Analysis} \textbf{11} (2012), 557-586. \vsm
\bibitem{appeul} C. Morosi, L. Pizzocchero, \textsl{On approximate solutions for
the Euler and Navier-Stokes equations}, Nonlinear Analysis
\textbf{75} (2012), 2209-2235.  \vsm
\bibitem{cok} C. Morosi, L. Pizzocchero, \textsl{On the constants
in a basic inequality for the Euler and NS equations}, Appl. Math.
Lett. \textbf{26} (2013), 277-284.
\bibitem{apprey} C. Morosi, L. Pizzocchero, \textsl{On the Reynolds
number expansion for the Navier-Stokes equations},
Nonlinear Analysis \textbf{95} (2014), 156-174.
\bibitem{smo}
C. Morosi, L. Pizzocchero, \textsl{Smooth solutions of the Euler and
Navier-Stokes equations from the a posteriori
analysis of approximate solutions}, Nonlinear Analysis \textbf{113} (2015), 298-308.
\vsm
\bibitem{Rau}
G. Raugel, G.R. Sell,
\textsl{Navier-Stokes equations on thin 3D domains. I: global attractors and global regularity of
solutions}, Journal of the American Mathematical Society \textbf{6} (1993), 503-568.
\vsm
\bibitem{Rob} J.C. Robinson, W. Sadowski, \textsl{Numerical verification of
regularity in the three-dimensional Navier-Stokes equations for
bounded sets of initial data}, Asymptot. Anal. \textbf{59} (2008),
39-50.
\vsm
\bibitem{Rob2}
J. C. Robinson, W. Sadowski,
\textsl{The regularity problem for the three-dimensional Navier-Stokes equations}, in:
Partial Differential Equations and Fluid Mechanics,
London Math. Soc. Lecture Note Ser. \textbf{364}, 185-206, Cambridge Univ. Press (2009).
\vsm
\bibitem{RSS}
J. C. Robinson, W. Sadowski, R. P. Silva,
\textsl{Lower bounds on blow up solutions of the three-dimensional Navier–Stokes equations
in homogeneous Sobolev spaces}, J. Math. Phys. \textbf{53} (2012), 115618, 15pp.
\vsm
\bibitem{Tay} G.I. Taylor, A.E. Green,
\textsl{Mechanism of the production of small eddies from large
ones}, Proc. R. Soc. Lond. A \textbf{158} (1937), 499-521.
\vsm
\bibitem{Tem} R. Temam, \textsl{Local existence of $C^\infty$ solutions of the Euler equation
of incompressible perfect fluids}, in ``Turbulence and Navier Stokes equation'',
Proceedings of the Orsay Conference,
Lecture Notes in Mathematics \textbf{565} (1976), 184-193.
\bibitem{gmpy} The GMPY Collaboration, ``Multiprecision arithmetic for Python'', see
\phantom{xxxxxxxxxx}
http://code.google.com/p/gmpy. This software is a wrapper for
GMP Multiple Precision Arithmetic Library, see http://gmplib.org.
}
\end{thebibliography}
\end{document}